\documentclass{amsart}

\newtheorem{theorem}{Theorem}[section]
\newtheorem{lemma}[theorem]{Lemma}
\newtheorem{proposition}[theorem]{Proposition}

\theoremstyle{definition}

\newtheorem{remark}[theorem]{Remark}
\newtheorem*{acknowledgement}{Acknowledgement}

\theoremstyle{remark}

\usepackage{amscd,amssymb}

\newcommand\mylabel[1]{\label{#1}}

\newcommand{\ZZ}{\mathbb{Z}}
\newcommand{\QQ}{\mathbb{Q}}

\newcommand{\CC}{\mathbb{C}}

\newcommand{\PP}{\mathbb{P}}
\renewcommand{\AA}{\mathbb{A}}
\newcommand{\GG}{\mathbb{G}}

\newcommand  {\shM}     {\mathcal{M}}

\newcommand  {\shL}     {\mathcal{L}}

\newcommand  {\foB}     {\mathfrak{B}}

\newcommand  {\foE}     {\mathfrak{E}}

\newcommand  {\foU}     {\mathfrak{U}}

\newcommand  {\foX}     {\mathfrak{X}}
\newcommand  {\foY}     {\mathfrak{Y}}

\newcommand  {\ab}      {{\operatorname{ab}}}

\newcommand  {\Alb}     {\operatorname{Alb}}

\newcommand  {\Aut}     {\operatorname{Aut}}

\newcommand  {\Br}      {\operatorname{Br}}

\newcommand  {\codim}   {\operatorname{codim}}

\newcommand  {\Div}     {\operatorname{Div}}

\newcommand  {\Fr}      {\operatorname{Fr}}

\newcommand  {\Hom}     {\operatorname{Hom}}

\newcommand  {\ind}     {\operatorname{ind}}

\renewcommand  {\k}     {\kappa}

\newcommand  {\lra}     {\longrightarrow}

\newcommand  {\maxid}   {\mathfrak{m}}

\newcommand  {\NS}      {\operatorname{NS}}

\renewcommand{\O}       {\mathcal{O}}

\newcommand  {\Pic}     {\operatorname{Pic}}

\newcommand  {\pr}      {\operatorname{pr}}

\newcommand  {\qM}      {\overline{M}}
\newcommand  {\ra}      {\rightarrow}

\newcommand  {\rank}    {\operatorname{rank}}
\newcommand  {\red}     {{\operatorname{red}}}

\newcommand  {\sep}     {{\operatorname{sep}}}

\newcommand  {\sh}      {{\operatorname{sh}}}

\newcommand  {\Spec}    {\operatorname{Spec}}
\newcommand  {\Spf}     {\operatorname{Spf}}

\newcommand  {\trdeg}   {\operatorname{trdeg}}

\def\mydate{\number\day\space\ifcase\month \or January\or February\or March\or 
April\or May\or June\or July\or
August\or September\or October\or November\or December\fi \space\number\year}

\begin{document}

\title[Strong Franchetta Conjecture]
      {The strong Franchetta Conjecture in
       arbitrary characteristics}

\author[Stefan Schroer]{Stefan Schr\"oer}
\address{Mathematische Fakult\"at, Ruhr-Universit\"at, 
         44780 Bochum, Germany}
\curraddr{}
\email{s.schroeer@ruhr-uni-bochum.de}


\dedicatory{\mydate}

\begin{abstract}
Using Moriwaki's calculation of the $\QQ$-Picard group
for the moduli space of curves, I prove 
the strong Franchetta Conjecture in all characteristics. 
That is, the canonical class generates the group of rational
points on the Picard scheme for the generic curve of
genus $g\geq 3$. Similar results hold for
generic pointed curves. Moreover, I show that 
Hilbert's Irreducibility Theorem implies that
there are many other nonclosed points in the moduli space
of curves with such properties.
\end{abstract}

\maketitle
\tableofcontents

\section*{Introduction}
Let $M_g$ be the coarse moduli space of  smooth curves of genus
$g\geq 3$ over an arbitrary ground field $k$.
Deligne and Mumford \cite{Deligne; Mumford 1969}  showed that $M_g$
is an  irreducible algebraic scheme.
Let $\eta\in M_g$ its the generic point and $C=M_{g,1}\ra M_g$ 
the tautological curve.
The generic fiber $C_\eta$ is a smooth curve of genus $g$  over the
function field
$\kappa(\eta)$ of the moduli space
$M_g$. We call it the \emph{generic curve}.

Franchetta \cite{Franchetta 1954}  conjectured
$\Pic(C_\eta)=\ZZ K_{C_\eta}$.
Arbarello and Cornalba \cite{Arbarello; Cornalba 1987}  proved it 
over the complex numbers using Harer's   calculation \cite{Harer
1983} of the second homology for the mapping class group of Riemann
surfaces. The latter is a purely topological result. Later,
Arbarello and Cornalba 
\cite{Arbarello; Cornalba 1998} gave an algebro-geometric proof over the
complex numbers. Mestrano \cite{Mestrano 1987} and
Kouvidakis \cite{Kouvidakis 1991} deduced
the \emph{strong Franchetta Conjecture} over $\CC$, which
states that the rational points in the Picard scheme 
$\Pic_{C_\eta/\eta}$ are precisely the multiples of the canonical
class.

The first goal of this paper is to give an algebraic
proof for the  strong Franchetta
Conjecture in all characteristics $p\geq 0$.
The idea is to construct special stable curves
showing that any divisor class violating Franchetta's Conjecture
must be nontorsion.
Having this, we use Moriwaki's calculation \cite{Moriwaki 2001} of 
$\Pic(\overline{M}_{g,n+1})\otimes\QQ$ in
characteristic $p>0$ to infer the strong Franchetta Conjecture.
I also prove the Franchetta Conjecture for 
\emph{generic pointed curves}:
Their Picard groups are freely generated
by the canonical class and the marked points.
Actually, I use the pointed case as an essential step
in the proof for the unpointed case.

The second goal of this paper is to show that there are many other
nonclosed points $x\in M_{g,n}$ such that the marked points and the
canonical class generates $\Pic(C_x)$, at least up to torsion.
This seems to be new even in characteristic zero.
We shall see that over   uncountable ground fields, there
is an uncountable dense set of such points with
$\dim\overline{\left\{x\right\}}\leq 2$. This relies on 
\emph{Hilbert's Irreducibility Theorem} for function fields. The idea
is to view $C_\eta$ as the generic fiber of some fibered surface $Y$,
extend this to a family of fibered surfaces $\foY\ra S$, and
apply Hilbert's Irreducibility Theorem and the
\emph{Tate--Shioda Formula} to the resulting family
of  N\'eron--Severi scheme $s\mapsto \NS(\foY_s)$.
Such specialization arguments are problematic in characteristic
$p>0$, because ungeometric properties like   regularity behave
badly in families. However, we overcome these difficulties by using
the theory of \emph{geometric unibranch singularities}.

Here is a plan for the paper.
The first section contains some general facts on curves,
Picard schemes, and moduli spaces.
In Section 2 we examine 
curves of compact type and ordinary abelian
varieties. In Section 3 we prove that the Mordell--Weil group of
the generic Jacobian is torsion free.
This result is further improved in Section 4.
Section 5 contains the proof for   the strong Franchetta Conjecture.
As an application I deduce in Section 6 that the generic  curve in
characteristic $p=2$ does not admit a tamely ramified morphism to
the projective line. In Section 7 we construct an explicit stable
curve $X$ of genus $g$ over the rational function field $F$
so  that the $\Pic_{X/F}^0$-torsor
$\Pic^1_{X/F}$  has order $2g-2$
in the Weil--Ch\^atelet group.

The next two sections contain some general result on schemes:
In Section 8, we show that geometrically unibranch schemes
are good substitutes for normal schemes in characteristic $p>0$.
In Section 9 we use Hilbert's Irreducibility Theorem,
which comes from Galois theory, to study Picard numbers
in families of proper schemes.
We shall apply these results in Section 10: There we first
construct fibered surfaces with small Picard numbers over
transcendental extension fields. Specializing them, we
show that there are many nonclosed points   $x\in M_{g,n}$
such that $\Pic(C_x)$ has the same rank as $\Pic(C_\eta)$.
The last section contains a list of open problems.

\begin{acknowledgement}
This paper comprises my Habilitationsschrift at the
Ruhr-Universit\"at.
I whish to thank my teacher Hubert Flenner for 
many discussions, constant encouragement, and all the support
over the years.
This research originated in discussions with Jan Nagel.
I thank him for many stimulating discussions
and   suggestions. I also thank J\'anos Koll\'ar 
for helpful discussions. 
This research was in part carried out at the Alfr\'ed R\'enyi
Institute in Budapest and the Departement of Mathematics in Toronto.
I am grateful to  these institutes for their hospitality.
\end{acknowledgement}

\section{Preliminaries}

Let us collect some well-known facts on algebraic curves and their
Pi\-card schemes. 
Throughout this paper we fix an arbitrary ground field  $k$ of
characteristic $p\geq 0$.
We will, however, also deal with algebraic schemes over
extension fields $k\subset F$.
A \emph{curve} over $F$ is a proper 1-dimensional
$F$-scheme $X$ with $F=H^0(X,\O_X)$. Each curve comes along with its 
\emph{Picard scheme}
$\Pic_{X/F}$, which is a smooth group scheme of finite type.
We denote by $\Pic(X/F)=\Pic_{X/F}(F)$ its group of rational
points. 
If $X$ contains a rational point, this is nothing but the Picard
group of $X$. In general we have an exact sequence
$$
0\lra\Pic(X)\lra\Pic(X/F)\lra\Br(F)\lra\Br(X),
$$
where $\Br(F)$ is the Brauer group
(\cite{GB III}, Corollary 5.3). 
Tsen's Theorem implies the following
(compare \cite{GB III}, Theorem 1.1):

\begin{lemma}
\mylabel{Tsen}
Suppose that $k$ is algebraically closed and that
$F$ is the function field of a smooth algebraic curve over $k$.
Then $\Pic(X)=\Pic(X/F)$.
\end{lemma}

Note that we do not loose rational points
by passing to field extensions:

\begin{lemma}
\mylabel{extension}
Let $F\subset F'$ be a field extension, and $X'=X\otimes_F F'$  the
induced curve. Then the canonical map $\Pic(X/F)\ra\Pic(X'/F')$ is
injective.
\end{lemma}

\proof
We have $\Pic_{X'/F'}=\Pic_{X/F}\otimes F'$,
hence the fiber over the rational point  $0\in \Pic_{X/F}$ is
nothing but the rational point $0\in\Pic_{X'/F'}$.
\qed

\medskip
Let $\Pic^0_{X/F}\subset\Pic_{X/F}$ be the connected component
of the origin. Equivalently, $\Pic^0_{X/F}$ is 
the subgroup scheme given by numerically trivial line bundles.
We sometimes write $J=\Pic^0_{X/F}$ and call it the \emph{Jacobian}
of $X$. If
$X$ is geometrically irreducible, the invertible sheaves of degree
$d$ comprise the other connected component
$\Pic^d_{X/F}\subset\Pic_{X/F}$.  These are torsors under the group
scheme $\Pic^0_{X/F}$, in other words, elements of the
\emph{Weil--Ch\^atelet} group 
$H^1(k,\Pic^0_{X/F})$. Here cohomology is with respect
to the fppf topology. However, \'etale cohomology gives the
same result because $\Pic^0_{X/F}$ is smooth (\cite{GB III},
Theorem 11.7).  The torsor $\Pic^d_{X/F}$ is trivial if and only
if it contains a rational point.
If $\Pic^0_{X/F}$ is an abelian variety,  $\Pic^0(X/F)$ is also
called the \emph{Mordell--Weil group}.

Given  an integer $g\geq 3$,
let $M_g$ be the \emph{coarse moduli space} of  smooth curves of
genus $g$  over $k$. Mumford \cite{Mumford; Fogarty; Kirwan 1993} 
showed that this is an
algebraic $k$-scheme. Moreover, 
Deligne and Mumford \cite{Deligne; Mumford 1969} proved
that $M_g$ is irreducible. If $k$ is algebraically closed,
its closed points correspond to   isomorphism classes of smooth
curves over $k$. In any case, the coarse moduli space $C=M_{g,1}$
of pointed smooth curves defines a tautological  family of curves
$C\ra M_g$.

We call $C_\eta$ the \emph{generic curve of genus $g$},  where
$\eta\in M_g$  is the generic point. 
This is justified as follows:
Let $\bar{F}$ be an algebraic closure  of the function field
$F=\k(\eta)$, and
$\bar{X}$ a  smooth curve of genus $g$ over $\bar{F}$ 
corresponding to  the geometric point $\Spec(\bar{K})\ra M_g$. Then
$C_\eta\otimes_F\bar{F}\simeq\bar{X}/\Aut(\bar{X})$.
However, the generic curve of genus $g\geq 3$  has trivial
automorphism group (see \cite{Poonen 2000} for an algebraic proof).
So  the generic curve $C_\eta$ is indeed a smooth curve of genus 
$g$ over
the function field $F=k(\eta)$.  In contrast,
the generic curve in
genus two is the the projective line.

As often, it will be important to consider  \emph{stable pointed
curves} of genus $g\geq 3$ as well. Over an algebraically closed
field, stability means that $X$ has only ordinary double points,
and each smooth rational component contains either two double
points, or two marked points, or a double points and a marked
point. Such curves form a coarse moduli space
$\overline{M}_{g,n}$, which is irreducible and projective.
Here $g$ is the (arithmetic) genus, and $n$ is the  number of marked
points
$x_i\in X$.
We also have a tautological family
$\pi:\overline{M}_{g,n+1}\ra\overline{M}_{g,n}$ 
sending an $(n+1)$-pointed stable curve
$(X,x_1,\ldots,x_{n+1})$ to the $n$-pointed stable curve 
$(X',x'_1,\ldots,x'_{n})$ obtained by forgetting the last
marked point and contracting the possibly occurring
rational component in
$(X,x_1,\ldots,x_n)$ violating stability. 

Let $\eta_n\in\overline{M}_{g,n}$ 
be the generic point and $C=\overline{M}_{g,n+1}$  the tautological
family of curves. 
The generic fiber $C_{\eta_n}$ is a smooth curve of
genus
$g$, which is endowed with marked rational points $c_1,\ldots,c_n\in
C_{\eta_n}$. We call $C_{\eta_n}$ the \emph{generic $n$-pointed
curve of genus $g$}.

\section{Curves of compact type}

Let $k\subset F$ be a field extension. A stable curve
$X$ over $F$ is called of \emph{compact type} if $\Pic^0_{X/F}$ is 
an abelian variety. 
Equivalently, the map $\Pic^0_{X/F}\ra\Pic^0_{\tilde{X}/F}$ 
induced by
the normalization $\tilde{X}\ra X$ is injective.
If $F$ is algebraically closed,
the condition means that the irreducible components $X_i\subset X$
are smooth, and  the configuration of the $X_i$ is tree-like
(\cite{Bosch; Luetkebohmert; Raynaud 1990}, Chapter 9, Proposition
10). Such curves have nice properties:

\begin{lemma}
\mylabel{transcendental}
Let $F\subset F'$ be a purely transcendental field extension,
$X$ a stable curve over $F$, and $X'=X\otimes F'$ the induced curve.
If $X$ is of compact type, then the preimage map 
$\Pic(X)\ra\Pic(X')$ is bijective.
\end{lemma}

\proof
By the usual limit argument (\cite{EGA IVd}, Theorem 8.5.2), we may
assume that $F'$ is finitely generated. Applying induction,
we reduce to the case $\trdeg(F')=1$. Now we may view $F'$
as the function field of $\PP^1_F$.

The task is to check that $\Pic(X)\ra\Pic(X')$ is surjective.
First we do the special case that $X$ is smooth.
Then $X\times\PP^1_F$ is factorial, so the restriction map 
$\Pic(X\times\PP^1_F)\ra\Pic(X')$ is surjective.
Using $\Pic(X\times\PP^1_F)=\Pic(X)\oplus\ZZ$ 
we deduce that $\Pic(X)\ra\Pic(X')$ is surjective as well.

Now suppose $X$ is arbitrary, and let $\tilde{X}\ra X$ be the
normalization. The maps $\Pic(X)\ra\Pic(\tilde{X})$ and
$\Pic(X')\ra\Pic(\tilde{X}')$ are  bijective, because $X$
is of compact type. Applying the preceding special case to $\tilde{X}$,
we infer that $\Pic(X)\ra\Pic(X')$ is surjective.
\qed

\medskip
Recall that, in characteristic $p>0$,
an abelian variety $A$ over $F$ has a \emph{$p$-rank} $f$
defined by $\Hom(\mu_p,A\otimes\overline{F})=(\ZZ/p\ZZ)^f$. 
The abelian variety is
called \emph{ordinary} if $f=\dim(A)$.
This condition means that the group of geometric $p$-torsion points
in $A(\overline{F})$ is as large as possible, namely 
$(\ZZ/p\ZZ)^{\dim(A)}$.
In characteristic zero,   abelian varieties are ordinary by
definition. The following specialization result will play a role
in the sequel:

\begin{lemma}
\mylabel{generization}
Suppose $B$ is a discrete valuation ring with residue field
$F$ and field of fractions $Q$. Let $Z\ra\Spec(B)$ be a 
relative stable curve
such that $Z_0$ is of compact type and $\Pic^0_{Z_0/F}$ is ordinary.
If $\Pic^0(Z_0/F)$ is torsion free, then $\Pic^0(Z_\eta/Q)$ 
is torsion free as well.
\end{lemma}

\proof
Consider the relative generalized Jacobian $J=\Pic^0_{Z/B}$.
This is a separated group scheme of finite type over $B$
by \cite{Deligne 1985}, Proposition 3.4. It is an abelian scheme
because $Z_0$ is of compact type.
Let $J[m]$ be the
relative kernel of the multiplication map $[m]:J\ra J$.
Then $J[m]$ is a finite flat group scheme over $B$ whose fibers have
length $m^{2g}$, where $g$ is the genus of $Z_0$.

We claim that
$J[m](Q)=0$. Indeed, suppose we have a  point
$x\in J[m](Q)$. 
We now assume that the characteristic $p>0$ is positive, the case
$p=0$ being similar.
Decompose $m=lq$, where $l$ is prime to $p$, and $q$ is a  power of
$p$.  By construction,
the geometric closed fiber of $J[m]$ is isomorphic to
$(\ZZ/l\ZZ)^{2g}\oplus(\ZZ/q\ZZ)^g\oplus\mu_{q}^g$,
where $\mu_{q}$ is the local group scheme of $q$-th roots of unity
(here we use that $J_0$ is ordinary).
The same holds for all
geometric fibers, because the number of  geometric connected
components is lower semicontinous in proper flat families
(\cite{EGA IVc}, Proposition 15.5.9). By
\cite{EGA IVd}, Proposition 15.5.1, each of the
$l^{2g}q^{g}$  geometric points 
in the closed fiber extend to disjoint sections
over the strict henselization
$B\subset B^{\sh}$, which defines $l^{2g}q^{g}$ 
different geometric points over the
generic fiber. It now follows
from \cite{SGA 1}, Expos\'e VIII, Theorem 4.1 that 
the closure   $\overline{\left\{x\right\}}\subset J[m]$ is disjoint
from the zero section, hence the specialization map
$J[m](Q)\ra J[m](F)$ is injective.
Since $J(F)$ is torsion free we have $J[m](Q)=0$.
\qed

\section{Torsion points}

Fix an integer $n\geq 0$ and a genus $g\geq 3$.
Let $C_{\eta_n}$ be the generic $n$-pointed
curve of genus $g$.
In this section we take care of torsion points:

\begin{proposition}
\mylabel{torsion free}
The group $\Pic^0(C_{\eta_n}/\eta_n)$ is torsion free.
\end{proposition}

The proof requires some preparation. 
First note that by Lemma \ref{extension}, we may 
replace the ground field $k$ by any extension field.
For the rest of the section, we assume that $k$  is algebraically
closed, and write $F=k(T)$ for the rational function field in one
indeterminate. 

\begin{proposition}
\mylabel{elliptic}
There is an  ordinary elliptic curve $E$ over  the rational function field $F$
satisfying $\Pic^0(E)=0$.
\end{proposition}

\proof
The idea is to use  special Halphen pencils.
Fix an ordinary elliptic curve $E_0$ over $k$.
Such a curve exists because there are only finitely many supersingular
elliptic curves (\cite{Hartshorne 1977}, Chapter IV, Corollary 4.23).
Let $x\in E_0$ be the origin and consider the closed embedding
$E_0\subset \PP^2_k$ defined via $\O_{E_0}(3x)\simeq\O_{E_0}(1)$.

Lines and quadrics
in $\PP^2_k$ are uniquely determined by their intersection with 
the cubic $E_0$, because the
restriction maps 
$H^0(\PP^2,\O_{\PP^2}(n))\ra H^0(E_0,\O_{E_0}(3nx))$ 
are bijective for $n=1,2$. 
Let $L\subset\PP^2_k$ be the unique line with $L\cap E_0=3x$.
Setting $E_\infty=3L$, we have $E_0\cap E_\infty=9x$. 
Consider the pencil of cubics
$E_t$, $t\in\PP^1_k$ generated by $E_0$ and $E_\infty$.
Then $E_t$ is integral for $t\neq\infty$. Indeed, any quadric $Q\subset E_t$
satisfies $Q\cap E_0=6x$, hence $Q=2L$. Moreover, any line $H\subset E_t$ 
has $H\cap E_0=3x$, hence $H=L$.

Now let $g':Y'\ra \PP^2_k$ be the blowing-up of the nonreduced center
$E_0\cap E_\infty$, and $f':Y'\ra\PP^1_k$ 
the induced fibration. The Jacobian of the generic fiber
$Y'_\eta$ is the desired elliptic curve $E$.
To see this, note that 
$R'={g'}^{-1}(E_0\cap E_\infty)$  is isomorphic 
to the projective line over the Artin algebra
$k[v]/(v^9)$. 
A local computation as in \cite{Morrison 1985}, page 417  reveals
that
$Y'$ is smooth except for a rational double point   $y'\in Y'$ of type $A_8$ 
lying on $R'$. Let 
$Y\ra Y'$ be the minimal resolution of this singularity.  Its
exceptional locus is a string of eight  smooth rational
$(-2)$-curves. If follows that the strict transform  
$R\subset Y$ of $R'_\red$ is
a smooth rational $(-1)$-curve.

Let $f:Y\ra \PP^1_k$ be the induced fibration and consider
its generic fiber $Y_\eta$.
We now check that $R$ generates $\Pic(Y_\eta)$. Indeed, the restriction
map  $\Pic(Y)\ra\Pic(Y_\eta)$ is surjective, and $\Pic(Y)$ 
is generated by the exceptional curves
for the sequence of blowing ups $g:Y\ra \PP^2_k$
and $g^*(L)$. The strict transform of $L$ is disjoint from $Y_\eta$,
because the strict transform of $3L$ is disjoint from $Y_\eta$, 
and we conclude that the restriction of $f^*(L)$
to $Y_\eta$ is a multiple of $R$.

Summing up, the Jacobian $E=\Pic^0_{Y_\eta/\eta}$ 
contains no rational point but
the origin. Then $\Pic^0(E)=0$, because $E$ is isomorphic to
its own Jacobian. 
To see that the elliptic curve
$E$ is ordinary, look at
the relative Jacobian  of $f:Y\ra\PP^1_k$ near $0\in\PP^1_k$.
Its closed fiber $E_0$ is ordinary. Since  this is an open
condition, the generic fiber $E$ is ordinary as well.
\qed

\begin{remark}
\mylabel{fiber}
(i) The rational double point $y'\in Y'$ maps
to $\infty\in\PP^1_k$. To see this,
write $\O_{\PP^2,x}^\wedge=k[[u,v]]$ so that $E_0,L\subset\PP^2_k$ 
correspond to $u=0$ and $v^3=u$, respectively.
Then $\O_{Y',y'}^\wedge=k[[u,v,w]]/(uw=v^9)$.
You easily check that the equation $v^3=u$  remains 
indecomposable inside $k[[u,v,w]]/(uw=v^9)$.
This equation defines the preimage ${g'}^*(L)\subset Y'$, which
decomposes as a Weil divisor into two irreducible components.
So the strict transform of $L$ in not Cartier,
hence passes through the singular point $y'\in Y'$.

(ii) The degenerate fiber $X_\infty\subset X$ is of type $\text{II}^*$
in Kodaira's notation, that is, it
corresponds to the root lattice $\widetilde{E}_8$.
To see this, note that the  string of eight
$(-2)$-curves in $X_\infty$ hits the strict transform of $L$
in precisely one point. Moreover, the strict transform
is not a $(-1)$-curve, because the
intersection form on $X_\infty$ is negative semidefinite.
Glancing at  Kodaira's classification 
(\cite{Kodaira 1963}, Theorem 6.2), we deduce 
that $X_\infty$ must be of type $\text{II}^*$.
\end{remark}

Next, we construct a stable curve $X$ of genus $g$
over the rational function field $F=k(T)$ as follows: 
Let $E_1,\ldots,E_g$ be copies of the elliptic curve $E$ from
Proposition \ref{elliptic}, and choose  
rational points $p_1,\ldots,p_g\in\PP^1_F$. Let 
$X=E_1\cup\ldots\cup E_g\cup\PP^1_F$ be the  curve obtained by
identifying  the rational point
$p_i\in\PP^1_F$ with the origins $0\in E_i$ for
$i=1,\ldots,g$. This curve is stable because $g\geq 3$. 
Note that we may view it as an $n$-pointed curve, simply by
choosing rational points in the rational component.
The curve $X$ has the following properties:

\begin{proposition}
\mylabel{properties}
The curve $X$ is of compact type, 
the abelian variety $\Pic^0_{X/F}$ is ordinary, and
its Mordell--Weil group is
$\Pic^0(X)=0$.
\end{proposition}

\proof
The normalization $\tilde{X}\ra X$ is the disjoint union of $g$
copies of
$E$ and a projective line. 
The canonical map $\Pic_{X/F}\ra \Pic_{\tilde{X}/F}$  is an
isomorphism,  so
$\Pic^0_{X/F}=E\times\ldots\times E$, and the result follows.
\qed

\medskip
\emph{Proof of Proposition \ref{torsion free}:}
Let $X=E_1\cup\ldots\cup E_g\cup\PP^1_F$ 
be the stable curve of genus $g$ over the
rational function field $F=k(T)$ constructed above.
Pick $n$ rational points $x_i\in X$  contained in the rational
component $\PP^1_F$ and disjoint from
the elliptic components $E_i$.
Let $\foY\ra\Spf(A)$ be the formal  versal deformation of the
$n$-pointed stable curve $X$. Then $A$ is a formal power series
ring  in $3g-3+n$ variables with coefficients in $F$.
Since $H^2(X,\O_X)=0$, we may extend any  ample invertible
$\O_X$-module to an   invertible $\O_{\foY}$-module. By
Grothendieck's Algebraization Theorem (\cite{EGA IIIa}, Theorem
5.4.5), the formal scheme $\foY$ is the formal completion of a
relative curve $Y\ra\Spec(A)$. 

Blowing up the closed point in $\Spec(A)$ and localizing at
the generic point of the exceptional divisor, we obtain
a discrete valuation ring $B$ dominating $A$.
Its field of fractions $B\subset Q$ is also the field of fractions
for $A$. The residue field  $B_0=B/\maxid_B$ is a purely
transcendental field extension of $F$.
Let $Z\ra\Spec(B)$ be the induced family of stable curves.

The classifying map $\Spec(B)\ra\overline{M}_{g,n}$ is dominant.
Hence $Z_Q=C_{\eta_n}\otimes Q$, which induces an  
injection $\Pic(C_{\eta_n})\subset\Pic(Z_Q)$.
We have $\Pic^0(Z_0)=0$ by construction, 
so Lemma \ref{generization} applies, and we conclude
that the group $\Pic^0(Z_Q)$ is torsion free.
\qed

\section{Another stable curve}
\mylabel{another}

We keep the notation from the previous section, such
that $C_{\eta_n}$ is the generic $n$-pointed curve of genus $g$.
Let $c_1,\ldots,c_n\in C_{\eta_n}$ be the marked points, and consider
the free abelian group 
$P=\ZZ c_1\oplus\ldots\oplus\ZZ c_n\oplus\ZZ K_{C_{\eta_n}}$.
The following is a key step in proving  Franchetta's Conjecture:

\begin{proposition}
\mylabel{cokernel}
The cokernel for the map $P\ra\Pic(C_{\eta_n})$ is torsion free.
\end{proposition}

This again depends on the existence of certain elliptic curves
over function fields. Let $F=k(T)$ be the rational function field,
and $E$ the ordinary elliptic curve with
$\Pic^0(E)=0$ from Proposition \ref{properties}.

\begin{proposition}
\mylabel{higher rank}
For each $m\geq 0$ there is a field extension $F\subset F'$
such that $E'=E\otimes F'$ satisfies $\Pic^0(E')=\ZZ^{\oplus t}$ 
for some $t\geq m$.
\end{proposition}

\proof
The case $m=0$ is trivial. We proceed by induction on $m$.
Suppose we already have a field extension
$F'$ such that $\Pic^0(E')=\ZZ^{\oplus t}$ 
with $t\geq m$.
Choose a basis $p_1,\ldots,p_t\in \Pic^0(E')$. Let $x\in E'$ 
be the origin, and $x_i\in E'$ the
rational points with $x-x_i\sim p_i$.

Now let $A$ be the henselization of the localized polynomial
algebra $F'[U]_{(U)}$, and consider the trivial family $Y=E'\otimes A$.
Let $S, S_i\subset Y$ be the sections corresponding to the rational points
$x,x_i\in E'$.
Let $u,v\in\O_{Y,x}$ be the regular parameter system corresponding
to $Y_0,S\subset Y$. Replacing $v$ by a different parameter $v'$, we obtain,
locally around $x$, a curve $S'\neq S$ with $S'\cap Y_0=\left\{x\right\}$.
According to \cite{EGA IVd}, Theorem 18.5.11, this locally defined
curve $S'$ defines a section $S_{t+1}\subset Y$ 
passing through
$x$ with $S_{t+1}\neq S$. Then the difference $S-S_{t+1}$ 
defines  a nonzero point
$p_{t+1}\in\Pic^0(Y_\eta)$ specializing to zero in $\Pic^0(Y_0)$. 
It follows that any nonzero multiple of $p_{t+1}$ lies outside  
the span of the $p_1,\ldots,p_t$.

Next, choose a subalgebra $B\subset A$ that is a localization of
a finite $F'[U]_{(U)}$-algebra at some $F'$-valued \'etale
point, so that the section
$S_{t+1}$ is defined over $B$.
Let $F''$ be the function field of $B$ and set $E''=E\otimes F''$.
Then $\Pic^0(E'')$ is torsion free according to
Proposition \ref{generization}.
By construction, it contains a free group of rank $t+1$. 

It remains to check that $\Pic^0(E'')$ is finitely generated.
Let $C$ be the normal curve over $F'$ corresponding to the
function field $F''$, and consider the regular proper surface $X=E'\times C$.
Then $\Pic^0_{X/F'}=\Pic^0_{E'/F'}\times\Pic^0_{C/F'}$, 
and the N\'eron--Severi group $\NS(X)=\Pic_{X/F'}/\Pic^0_{X/F'}$ 
is finitely generated.
Now view $E''$ as the generic fiber for the projection $X\ra C$.
The cokernel $\Pic^0(E'')/\Pic^0(E')$ 
injects into $\NS(X)$, and we conclude that
$\Pic^0(E'')$ is finitely generated.
\qed

\begin{remark}
The proof shows that we may choose $F\subset F'$ as a finitely generated
separable  field extension of transcendence degree $\leq m$.
\end{remark}

\medskip
We now use such elliptic curves to construct a geometrically
integral  $n$-pointed stable curve
of genus $g$.
Choose a field extension $F\subset F'$ so that $\Pic^0(E')=\ZZ^{\oplus t}$ 
for some $t\geq 2g-2+n$.
Pick $2g-2+n$ rational points 
$$
p_1,p_1',\ldots,p_{g-1},p_{g-1}',x_1,\ldots,x_n\in E'
$$ 
which are 
part of a basis for $\Pic(E')$.
Let $X$ be the stable curve of genus $g$ obtained by identifying
the pairs $p_i,p_i'\in E'$ for $i=1,\ldots,g-1$.
The rational points $x_i\in E'$ define $n$ rational points $x_i\in X$, which
we denote by the same letter.
The canonical morphism $\nu:E'\ra X$ is the normalization map for $X$.

\begin{proposition}
\mylabel{pullback}
We have $\nu^*(K_X)=\sum_{j=1}^{g-1}(p_j+p_j')$.
\end{proposition}

\proof 
This follows from duality theory for the finite
morphism $\nu:E'\ra X$.
See, for example, \cite{Reid 1994}, Proposition 2.3.
\qed

\medskip
\emph{Proof of Proposition \ref{cokernel}:}
Let $X$ be the $n$-pointed stable curve of genus $g$ over the
function field $F'$ constructed above.
As in the proof of Proposition \ref{torsion free}, we construct 
a discrete valuation ring $B$ and a stable curve $Z\ra\Spec(B)$ 
with regular total space such
that the following holds:
(i) The residue field $B_0=B/\maxid_B$ is a purely transcendental field
extension of $F'$, and the closed fiber is $Z_0=X\otimes B_0$.
(ii) If $B\subset Q$ denotes the function field, then
the classifying map $\Spec(Q)\ra\qM_{g,n}$ is dominant.

Now suppose we have a class $L\in\Pic(C_{\eta_n})$ such that 
$mL\in P$ for some $m\neq 0$.
Write $mL=\lambda K_{C_{\eta_n}}+\sum\lambda_i c_i$ 
for certain coefficients $\lambda,\lambda_i\in\ZZ$.
The task is to prove $L\in P$.
Using    the field extension $\k(\eta_n)\subset Q$,
we obtain a 
class $L_Q\in\Pic(Z_Q)$. It extends to a divisor $D\in\Div(Z)$, and we have
$$
mD = \lambda K_{Z/B}+\sum_{i=1}^{n}\lambda_i C_i 
$$
in $\Pic(Z)$, where $C_i\subset Z$ are the marked sections. 
This is because the divisor $Z_0$, being an integral fiber, supports
only principal divisors.
Pulling back to the normalization $\tilde{Z}_0=E\otimes B_0$
of the closed fiber $Z_0$,
we obtain 
$$
mD|_{\tilde{Z}_0} =
\lambda\sum_{j=1}^{g-1}(p_j+p_j')+\sum_{i=1}^{n}\lambda_i x_i 
$$
in $\Pic(E\otimes B_0)$.
Since the $p_j+p_j',x_i$ are part of a basis for $\Pic(E\otimes B_0)$
by Lemma \ref{transcendental}, all coefficients
$\lambda,\lambda_i$ are multiples of $m$.
Replacing $L$ by 
$L-\frac{\lambda}{m} K_{C_{\eta_n}}-\sum\frac{\lambda_i}{m} c_i$, 
we reduce to the
case that $mL=0$.
Then $L=0$ by Proposition \ref{torsion free}, and in particular  $L\in P$.
\qed

\section{The strong Franchetta Conjecture}

We come to the first main result of this paper:

\begin{theorem}
\mylabel{Franchetta}
Let $k$ be a field, and  $g\geq 3$ and $n\geq 0$ be integers.
Let $\eta_n\in \overline{M}_{g,n}$ the generic point in the moduli
space of $n$-pointed stable curves of genus $g\geq 3$, 
and $C=\overline{M}_{g,n+1}$ the tautological curve.
Then the  marked points $c_1,\ldots,c_n\in C_{\eta_n}$ and 
the canonical class $K_{C_{\eta_n}}$ 
freely generate $\Pic(C_{\eta_n}/\eta_n)$.
\end{theorem}

Before we prove this, let us recall the definition of certain
\emph{tautological classes} in $\Pic(\overline{M}_{g,n})\otimes\QQ$.
Let $\pi:\overline{M}_{g,n+1}\ra\overline{M}_{g,n}$ be the 
projection sending an $(n+1)$-pointed stable curve
$(X,x_1,\ldots,x_{n+1})$ to the $n$-pointed stable curve 
$(X',x'_1,\ldots,x'_{n})$ obtained by contracting any component in
$(X,x_1,\ldots,x_n)$ violating stability. The \emph{Hodge class}
$\lambda\in\Pic(\overline{M}_{g,n})\otimes\QQ$ is 
defined as the determinant of 
$\pi_*(\omega_\pi)$, where 
$\omega_\pi=\omega_{\overline{M}_{g,n+1}/\overline{M}_{g,n}}$ is the relative
dualizing sheaf.

We also have canonical sections
$s_i:\overline{M}_{g,n}\ra\overline{M}_{g,n+1}$ for $i=1,\ldots,n$
as follows.
These morphisms are best described on geometric points:
The section
$s_i$ sends an $n$-pointed stable curve $(X',x_1',\ldots,x_n')$ to
the
$(n+1)$-pointed stable curve
$(X,x_1,\ldots,x_{n+1})$ defined as follows:
We have $X=X'\cup\PP^1$, where the point $x_i'\in X'$
is identified with $\infty\in\PP^1$. The marked points
are 
$x_i=0\in\PP^1$, $x_{n+1}=1\in\PP^1$, and $x_j=x_j'$ for $j\neq i$.
The \emph{Witten classes} 
$\psi_i\in\Pic(\overline{M}_{g,n})\otimes\QQ$  are defined
as $\psi_i=s_i^*(\omega_\pi)$.

There are also \emph{boundary classes} 
$\delta_v\in\Pic(\overline{M}_{g,n})\otimes\QQ$, 
which are effective Weil divisors
supported on $\overline{M}_{g,n}-M_{g,n}$. They correspond to  
various  topological types of degeneration.
The idea now is to restrict tautological classes on $\qM_{g,n+1}$
to the generic curve $C_{\eta_n}\subset\qM_{g,n+1}$.

\begin{proposition}
\mylabel{restriction}
The subgroup $P\subset\Pic(C_{\eta_n})$ generated by the marked
sections and the canonical class contains the restriction 
of the tautological classes.
\end{proposition}

\proof
First note that geometric points $c$ in the generic $n$-pointed
curve $C_{\eta_n}$ correspond to $(n+1)$-pointed stable curves
$(X,x_1,\ldots,x_{n+1})$. The curve $X$ is smooth if $c$ is not a
marked point. On the other hand, if $c=c_i$ is a marked point, then
$X=X'\cup\PP^1$ with $x_i,x_{n+1}\in\PP^1$.
It now follows immediately from their definitions
in \cite{Moriwaki 2001}, Section 1 that the restrictions
$\delta_v|_{C_{\eta_n}}$ of the boundary classes are 
supported by the marked points $c_i\in C_{\eta_n}$.

How do   Hodge classes and Witten classes restrict to the generic
pointed curve? To see this, note that the fiber product
$C_{\eta_n}\times_{\overline{M}_{g,n+1}}\overline{M}_{g,n+2}$ 
is isomorphic to the blowing up of
$C_{\eta_n}\times_{\eta_n}C_{\eta_n}$  with respect to the centers
$(c_i,c_i)$ for
$i=1,\ldots,n$. This also follows from
the modular interpretation of geometric points
in $\overline{M}_{g,n+1}$. Indeed, the marked points
$c_i\in C_{\eta_n}$, $i=1,\ldots,n$ correspond to reducible stable
curves of the form
$(X'\cup\PP^1,x_1,\ldots,x_{n+1})$ with $x_i,x_{n+1}\in\PP^1$, and
the exceptional curve in the blowing up of 
$(c_i,c_i)\in C_{\eta_n}\times_{\eta_n}C_{\eta_n}$ is given by
$(n+2)$-pointed stable curves   $(X'\cup\PP^1,x_1,\ldots,x_{n+2})$
with $x_i,x_{n+1},x_{n+2}\in\PP^1$.
As a consequence we have $\lambda|_{C_{\eta_n}}=0$, 
and $\psi_i|_{C_{\eta_n}}=c_i$ for $i=1,\ldots,n$, and 
$\psi_{n+1}|_{C_{\eta_n}}=K_{C_{\eta_n}}$.
\qed

\medskip
\emph{Proof of Theorem \ref{Franchetta}:}
Set $P=\ZZ c_1\oplus\ldots\oplus\ZZ c_n\oplus\ZZ K_{C_{\eta_n}}$.
The canonical 
map $P\ra\Pic(C_{\eta_n})$ is injective. Indeed, this is clear for $n=0$.
In case $n\neq 0$, we have $\Pic(C_{\eta_n})=\Pic(C_{\eta_n}/\eta_n)$, 
and it suffices to construct
an $n$-pointed stable curve $X$ of genus $g$ where the marked
points and the canonical class are linearly independent.
We constructed such a curve in Section \ref{another}.

The   task is to prove surjectivity.
Fix a   point $L\in\Pic^d(C_{\eta_n}/\eta_n)$.
Having only quotient singularities, the normal scheme
$\qM_{g,n+1}$ is $\QQ$-factorial. Hence
some multiple $mL$ with $m>0$
extends to a Cartier divisor class $D$ on $\overline{M}_{g,n+1}$.
According to \cite{Moriwaki 2001}, Theorem 5.1, we have  
$D=a\lambda + \sum_{i=1}^{n+1} b_i\psi_i +\sum_v d_v\delta_v$  for
certain integral coefficients, at least after replacing
$m$ by a multiple.
Recall that 
$\lambda$ is the Hodge class,  $\lambda_i$ are the Witten classes,
and  $\delta_v$ are the boundary classes.
Restricting to the generic curve and using
Proposition \ref{restriction}, we infer 
$mL\in P$.

We now distinguish three case.
First, suppose $n\geq 1$. Then 
$\Pic(C_{\eta_n})=\Pic(C_{\eta_n}/\eta_n)$, 
and Proposition \ref{cokernel}
implies $L\in P$.
Second, suppose $n=0$ and $d=0$. Then there are only the two
tautological classes $\lambda$ and $\psi_1$ besides
the boundary classes. The equation
$$
mL=(a\lambda + \sum_{i=1}^{n+1} b_i\psi_i +\sum_v d_v\delta_v)|_{C_\eta}
$$
boils down to 
$mL=b_1K_{C_{\eta_n}}$ and therefore $mL=0$.
Now Proposition \ref{torsion free} ensures $L=0$.

So only the case $n=0$ and $d\neq 0$ remains.
Now we argue as follows:
The $\Pic^0_{C_\eta/\eta}$-torsor $\Pic^{d(2g-2)}_{C_\eta/\eta}$ 
contains both $(2g-2)L$ and $dK_{C_\eta}$. These rational points
differ by a   point in $\Pic^0(C_\eta/\eta)=0$, so
$(2g-2)L=dK_{C_\eta}$.
In other word, $L$ is a rational multiple of $K_{C_\eta}$.
Now consider the affine surjection  $C_{\eta_1}\ra C_{\eta}$ 
defined on geometric points by $(X,x_1,x_2)\mapsto (X,x_2)$.
We already saw that the canonical class
$K_{C_{\eta_1}}\in\Pic(C_{\eta_1})$ is a primitive
element.
Since $\Pic(C_\eta/\eta)\ra\Pic(C_{\eta_1})$ is injective,
the  canonical class $K_{\eta}\in\Pic(C_\eta/\eta)$ is primitive as well,
and we infer that   $L$ is an integral multiple of
$K_{C_{\eta}}$.
\qed

\begin{remark}
\mylabel{points}
By definition, the generic curves $C_{\eta_n}$  with $n\geq 1$
contain a rational point, hence
$\Pic(C_{\eta_n})=\Pic(C_{\eta_n}/\eta_n)$.   So a priori  there is
no difference between   weak and   strong Franchetta Conjectures.
\end{remark}

\begin{remark}
\mylabel{background}
Moriwaki's calculation \cite{Moriwaki 2001} of
$\Pic(\qM_{g,n})\otimes\QQ$ depends on the existence of  certain
simply connected coverings of
$\qM_{g,n}$ due to Looijenga \cite{Loojienga 1994}, 
Pikaart and de Jong \cite{Pikaart; de Jong 1994}, and
Boggi and Pikaart \cite{Boggi; Pikaart 2000}.

To my knowledge, the group $\Pic(\qM_{g,n})$ itself has not been
calculate yet, neither in characteristic zero nor in positive 
characteristics.
\end{remark}

\section{Application to tame coverings}

Let me give an application of the Strong Franchetta Conjecture
to tame coverings.
Belyi's Theorem \cite{Belyi 1979}  states that a complex curve is
defined over a number field if and only if it admits a map to
$\PP^1_\CC$ with at most three branch points.
Fulton showed that any smooth curve $X$ over a separably closed
field of characteristic 
$p\neq 2$ admits a branched covering  $X\ra\PP^1$  with tame
ramification (\cite{Fulton 1969}, Proposition 8.1).
Sa{\"\i}di used this to prove an analog of Belyi's Theorem 
in odd characteristics: A smooth curve $X$ over an
algebraically closed field of characteristic $p>2$ is defined over
a finite field if and only if it is a tamely ramified covering of
$\PP^1$ with at most three branch points
(\cite{Saidi 1997}, Theorem 5.6).
It is unknown to what extent these facts hold true
in characteristic $p=2$.
We have the following negative result for the generic $n$-pointed
curve $C_{\eta_n}$ of genus $g\geq 3$:

\begin{theorem}
\mylabel{wild}
Suppose the ground field $k$ has characteristic $p=2$.
Then any finite separable morphism $f:C_{\eta_n}\ra\PP^1_{\eta_n}$ 
has wild ramification.
\end{theorem}

\proof
Suppose on the contrary that  
$f:C_{\eta_n}\ra\PP^1_{\eta_n}$ has only tame ramification.
This means that, over the algebraic closure,
every ramification point  has odd ramification index.

Clearly, this implies that for all ramification points
$x\in C_{\eta_n}$, say with image $y\in \PP^1_{\eta_n}$,
the field extension
$\k(y)\subset\k(x)$ is separable.
Moreover, the localization of the fiber $f^{-1}(y)$ at $x$ is of the
form
$dx$ for some odd integer $d\geq 2$.
In turn,
the relative canonical class
$K_{C_{\eta_n}/\PP^1_{\eta_n}}=\sum(d-1)x$ is
$2$-divisible. The same holds for the canonical class of
$\PP^1_{\eta_n}$. Consequently $K_{C_{\eta_n}}$ is $2$-divisible,
contradicting Theorem \ref{Franchetta}.
\qed

\begin{remark}
\mylabel{spin}
The proof actually shows that $f:C_{\eta_n}\ra X$ has wild
ramification if
$X$ is any curve whose canonical class is 2-divisible in
the Picard scheme. The choice of an $L\in\Pic(X/\eta_n)$
with
$2L=K_X$ is sometimes called a \emph{spin structure}
or \emph{theta characteristic}, compare 
\cite{Arbarello et al 1985}, Chapter VI, Appendix B.
\end{remark}

\begin{remark}
Consider the unique supersingular elliptic curve $E$ in 
characteristic 
$p=2$, which has Weierstrass equation $y^2+y=x^3$ and invariant
$j(E)=0$.
Its automorphism group is a semidirect product of $\ZZ/3\ZZ$
and the quaternion group $\left\{\pm 1,\pm i,\pm j,\pm k\right\}$.
Therefore, this curve 
even admits a Galois covering $f:E\ra\PP^1$ of degree
$d=3$, which has three ramification points, all of
ramification index $e=3$.
Amusingly, this is the only elliptic curve in characteristic $p=2$
for which the existence of a  tamely ramified function is a priori
clear.
\end{remark}

It seems to be open whether a given smooth curve $X$ 
in characteristic $p\geq 0$ has a morphism
$f:X\ra\PP^1$ whose ramification indices $e\geq 2$  are all odd.
For example, we may ask this for the generic elliptic curve $E$
over $\QQ(j)$, which has Weierstrass equation
$$y^2+xy=x^3-\frac{1728}{j-1728}x-\frac{1}{j-1728}$$ and invariant
$j(E)=j$.
We can say at least the following:

\begin{proposition}
\mylabel{odd indices}
There is a  finite field extension $\QQ(j)\subset F$
and a morphism $f:E\otimes F\ra\PP^1_{F}$ of degree $d=4$
with four ramification points, all with ramification index $e=3$.
\end{proposition}

\proof
Let $\foE\ra S$ be the universal elliptic curve over
the punctured $j$-line $S=\AA^1_\QQ-\left\{0,1728\right\}$.
Consider the subscheme
$U\subset \Hom(\foE,\PP^1_S)$  of the
relative $\Hom$-scheme that parametrizes morphisms
$\foE_s\ra \PP^1_s$ of degree $d=4$ with four ramification points,
all with index $e=3$.
It follows from 
\cite{Fried; Klassen; Kopeliovich 2001}, Theorem 1, that
the projection $\pr:U\ra S$ is dominant.  Then the generic fiber
$\pr^{-1}(\eta)\subset U$ contains a closed point $u$, 
and the residue
field $F=\k(u)$ is the desired finite field extension of $\QQ(j)$.
\qed

\section{A stable curve with maximal index}

Recall that the \emph{index} $\ind(X)\in\ZZ$ 
of a proper curve $X$ over an arbitrary field
$F$ is the positive generator for the image of the degree map 
$\Pic(X/F)\ra\ZZ$. If $X$ is geometrically integral, $\ind(X)$ is
also the 
order of $\Pic^1_{X/F}$  in the Weil--Ch\^atelet group
$H^1(k,\Pic^0_{X/F})$.

The index of a stable curve of genus $g$ is a divisor of $2g-2$,
because there is always the canonical class.
It follows from Theorem \ref{Franchetta} that the
generic curve has maximal index $\ind(C_\eta)=2g-2$.
It would be interesting to construct special curves
with maximal index over fields of smaller transcendence degree. 
I do not know how to achieve this with smooth curves. The goal
of this section is to produce an explicit stable curve
with maximal index.

Throughout, we assume that our ground field $k$ is algebraically
closed, and let $F=k(T)$ be the rational function field.
Note that, in our situation, we have $\Pic(X/F)=\Pic(X)$ 
by Lemma \ref{Tsen}.

\begin{proposition}
\mylabel{genus one}
For each integer $d>0$ there is a smooth curve $Y$ of genus one over the
rational function field
$F$ with $\ind(Y)=d$.
\end{proposition}

\proof
Let $F\subset F'$ be a  cyclic field extension of degree $d$, say
with Galois group $G\simeq\ZZ/d\ZZ$. Choose an
elliptic curve
$E$ over $F$ containing a rational point $x\in E$ of order $d$.
Let $T\subset \Aut(E)$ be the cyclic subgroup of order $d$ generated by the
translation $T_x(e)=e+x$.

Set $E'=E\otimes_F F'$. Then $g\in G$ acts on $\Aut(E'/F')$ via 
$\phi\mapsto g\phi g^{-1}$. This action fixes the subgroup
$T$ pointwise, because $x\in E$ is a rational point.
We get
$$
H^1(G,T)=\Hom(G,T)\simeq\ZZ/d\ZZ.
$$
The latter identification  depends on the choice of a generator
$g_0\in G$. The inclusion $T\subset\Aut(E'/F')$ gives an inclusion 
$H^1(G,T)\subset H^1(G,\Aut(E'/F'))$. As explained in
\cite{Serre 1972}, Chapter II, \S1, Proposition 5, the set
$H^1(G,\Aut(E'/F'))$ may be viewed as the set of isomorphism classes of
twisted forms $Y$ of $E$ whose preimage
$Y'=Y\otimes F'$ is isomorphic to $E'$. Indeed,
$Y$ is the quotient of $E'$ by the $G$-action
$g\circ z_g$, where $z_g\in Z^1(G,\Aut(E'/F'))$ 
is a cocycle representing a given cohomology class.

Now choose the generator $g_0\mapsto T_x$ of $H^1(G,T)$ 
and consider the corresponding twisted
form $Y$ of $E$. Then $\ind(Y)$  divides $d$. Indeed, the reduced
divisor 
$D\subset E$ comprising  the multiples of $x\in E$ defines a
$G$-invariant divisor of degree
$d$ in $E'$, hence an element in $\Pic^d(Y/F)$.

For the converse, suppose we have a class in $\Pic^l(Y/F)$. It
corresponds to an invertible
$\O_Y$-module
$\shL$ of degree  $l$ by Lemma \ref{Tsen}. 
The isomorphism class of the induced invertible
$\O_{Y'}$-module
$\shL'$  is invariant under  both $g_0$  and $g_0\circ T_x$, hence
also under the translation $T_x$.
But $T_x^*(\O_{E'}(e))\simeq\O_{E'}(e-x)$  for any rational point
$e\in E'$,  so $T_x^*(\shL')\simeq\shL'$ implies $lx=0$. In turn,  
$d$ divides $l$, and we conclude $\ind(Y)=d$.
\qed

\begin{remark}
\mylabel{cyclic}
The cyclic field extension used above
exist. Indeed, if $d$ is prime to the characteristic $p$, the cyclic
extensions of degree
$d$ correspond to   \'etale cohomology
$H^1(F,\ZZ/d\ZZ)=F^*/(F^*)^d$ via the Kummer sequence
$0\ra\ZZ/d\ZZ\ra\GG_m\ra\GG_m\ra1$. If $d=p^m$, they
correspond to   \'etale cohomology   $H^1(F,\ZZ/d\ZZ)=F/\wp(F)$
via the Artin--Schreier sequence $0\ra\ZZ/p^m\ZZ\ra W_m(F)\ra
W_m(F)\ra 0$, where $W_m(F)$ is the sheaf of Witt vectors of length
$m$. In general, decompose $d=lp^m$ with $l$ prime to $p$, and let
$F'$ be the subextension  inside the maximal abelian extension
$F\subset F^\ab$ generated by the linearly disjoint
cyclic extensions corresponding to $l$ and $p^m$
(\cite{A 4-7}, Chapter V, \S 10, no.\ 4, Corollary 2 of Theorem 1).
\end{remark}

Now we assemble the desired stable curve:

\begin{proposition}
\mylabel{index}
There is a geometrically  integral stable curve $X$ of genus $g$
over the rational function fields
$F$ with $\ind(X)=2g-2$.
\end{proposition}

\proof
Let $Y$ be a smooth curve of genus one with $\ind(Y)=2g-2$, as in
Proposition \ref{genus one}.
By Lemma \ref{Tsen}, there is an  invertible $\O_Y$-module $\shL$ of
degree
$2g-2$. The Riemann--Roch Theorem gives $h^0(Y,\shL)=2g-2$,
hence $\shL$ comes  from an effective divisor $D\subset Y$ of length
$2g-2$. By Bertini's Theorem, the very ample sheaf $\shL$ has a
global section so that the corresponding divisor  $D\subset Y$ is
smooth  (see \cite{Jouanolou 1983}, Corollary 6.11).
Such a  subscheme is necessarily of the form $D=\Spec(F')$ for some 
separable field
extension
$F\subset F'$ of degree $2g-2$ because $\ind(Y)=\deg(\shL)$.

Let $F\subset F''\subset F'$ be  a
subextension of degree
$g-1$, such that $F''\subset F'$ has degree two.
The cocartesian diagram 
$$
\begin{CD}
\Spec(F') @>>> Y\\
@VVV@VVV\\
\Spec(F'')@>>> X
\end{CD}
$$
defines an integral curve $X$. The   curve  $X$ has a single 
cuspidal singularity, which breaks up
into $g-1$ nodal singularities over the algebraic closure.   
Therefore
$X$ is stable and geometrically integral. The exact sequence 
$$
0\lra\O_X\lra F''\oplus\O_Y\lra F'\lra 0
$$ 
gives an exact sequence
$$
0\lra H^0(\O_X)\lra F''\oplus H^0(\O_Y)\lra F'\lra H^1(X,\O_X)\lra
H^1(Y,\O_Y)\lra 0,
$$
so $h^1(X,\O_X)=g$.  
Since $Y\ra X$ is birational,   the map $\Pic(X)\ra\Pic(Y)$ is
surjective, and its kernel consists of numerically trivial sheaves.
We conclude that $\ind(X)=\ind(Y)=2g-2$ holds.
\qed

\begin{remark}
Unfortunately, the curve $X$ is not of compact type and contains
no rational point.
It is therefore difficult to analyze specialization of points
in the Picard scheme of the versal deformation for $X$.
In particular, the indices of   curves occurring in the
versal deformation are hard to control.
\end{remark}

\section{Geometrically unibranch singularities}

Our next goal is to find, beside the generic point
$\eta_n\in M_{g,n}$, other nonclosed points $x\in M_{g,n}$ 
whose curve $C_x$ has the property that the marked section
and the canonical class generate $\Pic(C_x)$,  at least up to
torsion. This will occupy the remaining sections.
But first we have to remove a problem in characteristic $p>0$,
which occurs over and over again, namely:
Normality is not necessarily preserved under 
inseparable field extensions. An annoying consequence
is, for instance, that the Jacobian $\Pic^0_{X/F}$
of a proper normal scheme $X$ is not necessarily proper.

In this section we sidestep such  problems using
geometrically unibranch schemes instead of normal schemes.
This property is stable under field extensions and, as we shall see,
almost as good as normality. Recall that a local ring $A$ is
called
\emph{geometrically unibranch} if it is irreducible, and the
normalization $\tilde{A}$ of the reduction
$A_\red$ induces a bijection
$\Spec(\tilde{A})\ra\Spec(A)$  whose residue field
extension is purely inseparable.
Equivalently, the strict henselization $A\subset A^\sh$ is
irreducible
(\cite{EGA IVd}, Proposition 18.8.15). A scheme is geometrically
unibranch if all its local rings are geometrically unibranch. The
nice thing about such schemes is that we may change the structure
sheaf $\O_X$ as long as the topological space remains the same:

\begin{lemma}
\mylabel{homeomorphism}
\mylabel{homeomorphism and unibranch}
Let $f:X\ra Y$ be an integral universal homeomorphism of schemes.
Then $Y$ is geometrically unibranch if and only if $X$ is
geometrically unibranch.
\end{lemma}

\proof
To check this we may assume that $X=\Spec(B)$ and $Y=\Spec(A)$ are
local. We have $B^\sh=B\otimes_A A^\sh$ according to \cite{EGA IVd},
Proposition 18.8.10 and Remark 18.8.11.
Hence $A^\sh$ is irreducible if and only if $B^\sh$ is irreducible,
because $A\ra B$ is a universal homeomorphism.
In turn, $A$ is geometrically unibranch if and only if $B$ is
geometrically unibranch.
\qed

\medskip
Let us now consider Jacobians. Surely, Jacobians of geometrically
unibranch proper schemes may be nonproper. However, things are
not too bad:

\begin{proposition}
\mylabel{proper image}
Let $X$ be  proper $F$-scheme in characteristic $p>0$. If $X$ is
geometrically unibranch, then 
there is an integer
$n\geq 0$ such that the schematic image of the multiplication map
$[p^n]:\Pic^0_{X/F}\ra\Pic^0_{X/F}$ is proper.
\end{proposition}

\proof
Let $H_n\subset\Pic^0_{X/F}$ be the image of the
multiplication map $[p^n]$. Being the image of a homomorphism of
algebraic group schemes, the embedding
$H_n\subset\Pic^0_{X/F}$ is   closed   by \cite{SGA 3a},
Expos\'e VI$_B$, Proposition 1.2.

Note that  we may replace the ground field $F$ by its algebraic
closure $\bar{F}$.
Indeed, the scheme $\bar{X}=X\otimes\bar{F}$ is geometrically
unibranch by
\cite{EGA IVb}, Proposition 6.15.7. Moreover, the image $H_n$
of
the multiplication map
$[p^n]$ on $\Pic^0_{X/F}$ is proper if and only if 
the induced image $\bar{H}_n=H_n\otimes\bar{F}$ is proper,
according to \cite{SGA 1}, Expos\'e VIII, Corollary 4.8.

Next, let
$\tilde{X}\ra X$ be the normalization of $X_\red$.
The  $\tilde{X}$ is geometrically normal, hence
$\Pic^0_{\tilde{X}/F}$ is proper by \cite{FGA VI}, Theorem 2.1. 
Let $\tilde{X}\ra\tilde{X}^{(p^m)}$ be the $m$-fold relative
Frobenius map, which is a finite universal homeomorphism.
It admits a factorization
$$
\tilde{X}\lra X\lra\tilde{X}^{(p^m)}
$$
for some $m\geq 0$ by
\cite{Kollar 1997}, Proposition 6.6, because
$\tilde{X}\ra X$ is a finite universal homeomorphism.
Applying this result again, this time to the finite universal
homeomorphism
$X\ra\tilde{X}^{(p^m)}$, we
obtain a factorization
$$
X\lra\tilde{X}^{(p^m)}\lra X^{(p^n)}
$$
for some $n\geq 0$. 
Consequently, the image of $\Pic^0_{\tilde{X}^{(p^m)}/F}$ in
$\Pic^0_{X/F}$ contains the image of
$\Pic^0_{X^{(p^n)}/F}$, which is nothing but $H_n$.
Using that $\tilde{X}^{(p^m)}$ is   normal, we conclude
that $H_n$ is proper.
\qed

\medskip
The next result tells us that geometrically unibranch schemes
behave well in families:

\begin{proposition}
\mylabel{unibranch constructible}
Let $S$ be an integral noetherian scheme in characteristic $p>0$,
and $f:X\ra S$ a morphism of finite type
whose generic fiber $X_\eta$
is geometrically unibranch.
Then there is a nonempty open subset $U\subset S$ such that all
fibers $X_s$, $s\in U$ are geometrically unibranch.
\end{proposition}

\proof
Using
generic flatness (\cite{EGA IVc}, Theorem 11.1.1),
we may replace $S$ by an open subset and assume that $f$ is flat.
By the usual limit argument,
there is a finite morphism $X'\ra X$ such that
$X'_\eta\ra X_\eta$ coincides with the generic fiber of the
normalization for
$X_\red$. By assumption, 
$X'_\eta\ra X_\eta$ is a
universal homeomorphism, hence its geometric fibers are connected.
According to
\cite{EGA IVc}, Proposition 15.5.1, the whole map $X'\ra X$ has
geometrically connected fibers, hence is a
universal homeomorphism.
We may replace $X$ by $X'$ and assume that $X_\eta$ is normal.

Next, we define for each integer $n\geq 0$ a
scheme
$X_n$ as follows. The underlying topological space for $X_n$ is the
underlying topological space for $X$. If $W\subset X$ is an affine
open subset, we set
$$
\Gamma(W,\O_{X_n})=\Gamma(W,\O_X) \cap Q(L^{p^n}).
$$
Here $Q$ is the function field of $S$ and $L$ is the function field
of $X$. Then each $(X_n)_\eta$ is a normal scheme and the projection
morphisms
$X_n\ra S$ are proper. The descending sequence of fields 
$L\supset L^p\supset\ldots$ 
defines a sequence
of   universal homeomorphisms
$X_0\ra X_1\ra\ldots $ with $X_0=X$.

According to \cite{Huebl; Kunz 1994}, Theorem 1, there is an integer
$n\geq 0$ such that the geometric fiber $Y_\eta$ of $Y=X_n$ is
geometrically normal over $Q$.
By \cite{EGA IVc}, Theorem 12.2.4, the set $W\subset Y$ of 
geometrically normal points  is open.  According to Chevalley's
Theorem (\cite{EGA I}, Theorem 7.1.4), the image $f(W)\subset S$ is 
constructible. We have $\eta\in f(W)$ by construction.
Being constructible, the set $f(W)$ contains a nonempty
open subset $U\subset S$.
The fibers $Y_s$, $s\in U$ are geometrically
unibranch. Since $X_s\ra Y_s$ is a finite universal
homeomorphism,  we may apply Lemma \ref{homeomorphism} and conclude
that the $X_s$, $s\in U$
are geometrically unibranch as well.
\qed

\medskip
The other direction holds as well:

\begin{proposition}
\mylabel{nonunibranch constructible}
Let $S$ be an integral noetherian scheme,
and $f:X\ra S$ a morphism of finite type
whose generic fiber $X_\eta$
is not geometrically unibranch.
Then there is a nonempty open subset $U\subset S$ such that all
fibers $X_s$, $s\in U$ are not geometrically unibranch.
\end{proposition}

\proof
Clearly, we may assume that $S$ is affine and  that $X$ is reduced.
Like in the proof for Proposition \ref{unibranch constructible},
there
is a finite morphism $X'\ra X$ such that $X'_\eta\ra X_\eta$ is the
normalization. Shrinking $S$, we may assume that  the fibers $X'_s$
have no embedded component by
\cite{EGA IVc}, Proposition 9.9.2.
According to \cite{EGA IVc}, Corollary 9.7.9,
the set $Z\subset X$ over which the geometric fibers of $X'\ra X$ 
are not connected is constructible.
Replacing $X$ by some open subset, we may assume that $Z\subset X$ is
closed.

The set $x\in X$   where $\codim_x(Z_{f(x)},X_{f(x)})=0$ is
constructible by
\cite{EGA IVc},  Proposition 9.9.1, and disjoint
from $X_\eta$. By Chevalley's Theorem,
there is a nonempty open subset $U\subset S$ such that
$\codim_x(Z_{f(x)},X_{f(x)})>0$ holds for all 
$x\in f^{-1}(U)$. In other words, $Z_s$ contains no generic point from
$X_s$ for  
$s\in U$. 
We infer that the fibers $X_s$
are not geometrically unibranch for $s\in U$.
\qed

\begin{remark}
\mylabel{constructible}
In Grothendieck's terminology, the preceding two result tell us
that the property of being geometrically unibranch is
\emph{constructible}
(compare \cite{EGA IVc}, Definition 9.2.1). It then follows from
\cite{EGA IVc}, Corollary 9.2.4, that given any morphism 
$f:X\ra S$ of finite
presentation, the set $E\subset S$ of points whose fiber $X_t$ is
geometrically unibranch is locally constructible.
\end{remark}

\medskip
Now suppose $X$ be a proper $F$-scheme. 
Recall that the \emph{N\'eron--Severi group scheme} is the quotient
$\NS_{X/F}=\Pic_{X/F}/\Pic^0_{X/F}$. This is an \'etale group
scheme over $F$, whose points are given by finite separable field
extensions (\cite{SGA 3a}, Expos\'e VI$_A$, 5.5). Its group
of rational points is denoted by
$\NS(X/F)$, which is a finitely generated abelian group
(\cite{SGA 6}, Expos\'e XIII, Theorem 5.1). The
\emph{N\'eron--Severi} group is the subgroup
$\NS(X)$ generated by $\Pic(X)$.

\begin{lemma}
\mylabel{Neron-Severi finite}
The inclusion $\NS(X)\subset\NS(X/F)$ has finite index.
\end{lemma}

\proof
Fix a rational point $z\in\NS(X/F)$. The corresponding connected
component
$T\subset\Pic_{X/F}$ is a torsor under $J=\Pic^0_{X/F}$, and we have
to show that its order in $H^1(k,J)$ is finite. Indeed, then $T^n$
contains a rational point $x$ for some $n>0$, and $mx$ comes from an
invertible sheaf  for some $m>0$ because $\Br(F)$ is torsion.

Let $F\subset F'$ be some finite field extension such that
$T'=T\otimes F'$ acquires a rational point $y$.
Passing to some multiple, we may also assume that $y$
comes from some invertible sheaf $\shM$.

First, consider the special case that $F'$ is separable.
Then we may assume that $F'$ is Galois as well, say
with Galois group $G$. Then $\bigotimes_{g\in G} g^*\shM$ is a
$G$-invariant invertible sheaf whose rational point lies on
${T'}^n$, hence descends to the desired rational point $x\in T^n$.

Next, consider the special case that $F'$ is purely inseparable.
Then the projection $X'=X\otimes F'\ra X$ is a finite universal
homeomorphism of $F$-schemes, so we have a factorization
$
X'\ra X\ra {X'}^{(p^n)}
$
for some $n\geq 0$. It follows that ${\shM}^{\otimes p^n}$ comes
from the  invertible
$\O_X$-module $\shL=\shM^{(p^n)}|X'$.

The general case follows from the special cases, because
any finite field extension is given by a   purely
inseparable extension followed by a separable extension.
\qed

\begin{proposition}
\mylabel{finite}
Let $f:X\ra Y$ be a universal homeomorphism of proper $F$-schemes
in characteristic $p>0$.
Then kernel and cokernel of the induced map $\NS(Y/F)\ra\NS(X/F)$
are finite $p$-groups.
In particular, we have $\rho(Y)=\rho(X)$.
\end{proposition}

\proof
y \cite{Kollar 1997}, Proposition 6.6, we have
a factorization $X\ra Y\ra X^{(p^n)}$ for some 
$n\geq 0$. Observe that the cartesian diagram
$$
\begin{CD}
{X}^{(p^n)} @>>> X\\
@VVV @VVV\\
\Spec(F) @>>\Fr^n> \Spec(F).
\end{CD}
$$
induces a bijection
$\NS(X/F)\ra \NS(X^{(p^n)}/F)$. This is because
$\Fr^n$ is a universal homeomorphism, and
sections for the \'etale scheme $\NS_{X/F}$ correspond to
points $l\in\NS_{X/F}$ such that $\Spec\k(l)\ra\Spec(F)$ is a universal
homeomorphism (\cite{EGA IVd}, Corollary 17.9.3).
It follows that kernel and cokernel of
$\NS(X^{(p^n)}/F)\ra\NS(X/F)$ are
$p$-groups.
Hence the cokernel of $\NS(Y/F)\ra\NS(X/F)$ is a $p$-group.

Finally, suppose $L\in\NS(Y/F)$ vanishes on $X$. As above,  there
is a factorization $Y\ra X^{(p^n)}\ra Y^{(p^m)}$ for some $m\geq 0$.
Then $p^mL$ comes from a point $M\in\NS(X^{(p^n)}/F)$,
and $p^nM$ vanishes on $X$. It follows that some $p$-power 
annihilates $L$.
\qed

\begin{remark}
\mylabel{char zero}
The preceding result does not hold in characteristic zero.
To see an example, set $X=\PP^2_\CC$, such that $\Pic(X)=\ZZ$. 
According to 
\cite{Hartshorne 1977}, Exercise 5.9 on page 232, there is a first
order
infinitesimal extension
$X\subset Y$ with ideal $\O_X(-3)=\omega_X$ such that  $\Pic(Y)=0$.
\end{remark}

\section{Picard numbers in families and Hilbert sets}

Let $S$ be a scheme and $f:X\ra S$ a
proper morphism. How do the Picard numbers $\rho(X_s)$
vary in this family? The function $s\mapsto\rho(X_s)$ can be
very nasty. To
explore its nature we need the notion of Hilbert sets. 
Suppose  for simplicity that $S$ is integral and noetherian. Let $T$
be an integral scheme and
$T\ra S$  an
\'etale dominant morphism of finite type. The generic fiber $T_\eta$
is an integral scheme given by a finite separable field extension
of
$\k(\eta)$. We define $H_T\subset S$ as the set of all points
$s\in S$ where the fiber $T_s$ is integral and nonempty.
A subset $H\subset S$ is called a \emph{separable Hilbert set} if
there are finitely many $T_i\ra S$, $1\leq i\leq n$ as above
with 
$H=H_{T_1}\cap\ldots\cap H_{T_n}$. Clearly,  open sets
are separable Hilbert sets. Moreover, separable Hilbert subsets
are stable under generization.  

Usually, Hilbert sets are defined in terms of irreducible
polynomials. They occur in Galois theory in connection with
\emph{Hilbert's Irreducibility Theorem}. The latter states, in its
original form, that an irreducible polynomial
$P(T,X)\in\QQ[T,X]$ remains irreducible for infinitely many
rational specializations  $T=a/b$. Hilbert sets are studied in detail
in \cite{Fried; Jarden 1986}, Chapter 11 and 12.
Our definition follows Lang \cite{Lang 1983} and 
Serre \cite{Serre 1989}, which is better suited for algebraic
geometry.   For example, Hilbert sets are closely related to
N\'eron--Severi groups:

\begin{proposition}
\mylabel{countable}
Let $S$ be an integral noetherian scheme, $\eta\in S$ its generic
point, and $f:X\ra S$ a
proper morphism. Then there are countably many nonempty
separable Hilbert subsets
$H_i\subset S$ such that $\rho(X_s)=\rho(X_\eta)$ for all
$s\in\bigcap H_i$.
\end{proposition}

\proof
Replacing $S$ by some nonempty open subset, we may
assume that the relative Picard functor $\Pic_{X/S}$ is representable
by a scheme, according to \cite{SGA 6}, Expos\'e XII, Corollary 1.2.
After shrinking $S$ further, we may also assume
that there is a group scheme $J\ra S$ of finite type, together with
an open embedding $J\subset \Pic_{X/S}$, such that
$J_s=\Pic^0_{X_s/s}$ for all $s\in S$. This holds by 
\cite{FGA VI}, Lemma 1.2.
By generic flatness, we may also assume that $J\ra S$ is flat.

Consider the group scheme $\NS_{X/S}$. 
By Chow's Lemma, there is a projective surjective morphism
$X'\ra X$ such that the composite map $f':X'\ra S$ is projective.
Shrink $S$ so that $\NS_{X'/S}$ is representable by a group scheme.
Then $\NS_{X'/S}$ contains only countably many irreducible component.
Indeed, this follows  from the existence of the Hilbert
scheme (see \cite{Bosch; Luetkebohmert; Raynaud 1990}, Chapter 8,
Theorem 5).
Using that $\Pic_{X/S}\ra\Pic_{X'/S}$ is of finite type 
(\cite{SGA 6}, Expos\'e XIII, Theorem 3.5), we conclude that
$\NS_{X/S}$ has only countably many irreducible components as well.

Let $A_i\subset \NS_{X/S}$ be such an irreducible component.
If 
$A_i\ra S$ is not dominant, choose a nonempty open subset
$H_i\subset S$ disjoint from the image of $A_i$. If $A_i\ra S$ is
dominant, let
$U_i\subset S$ be an open subset over which 
$A_i$ becomes flat, hence \'etale, and let
$H_i\subset U_i$ be the corresponding Hilbert subset.

Then $H=\bigcap H_i$ is the desired countable intersection of Hilbert
subsets. Indeed, for each $s\in H$, the points $a\in\NS(X_s/s)$
correspond to the points
$b\in\NS(X_\eta/\eta)$ via some dominant irreducible component $A_i$.
Moreover, $a$ is rational if and only if $b$ is rational.
In other words, we have constructed a  bijective specialization map
$\NS(X_\eta/\eta)\ra\NS(X_s/s)$ and conclude
$\rho(X_s)=\rho(X_\eta)$.
\qed

\medskip
How do countable intersections of Hilbert sets look like?
There seems to be no general answer. However, over  certain
ground fields we can say more;

\begin{lemma}
\mylabel{classical Hilbert}
Let $k$ be an uncountable field, $k\subset F$ a   finitely
generated field extension of transcendence degree $\geq 1$, and
$S$ an integral $F$-scheme of finite type of dimension $n\geq 1$.
Then any countable intersection $\bigcap H_i$  of nonempty
separable  Hilbert subsets $H_i\subset S$ contains uncountably many
closed points.
\end{lemma}

\proof
We may assume that $S=\Spec(R)$ is affine.
Let $R\subset Q$ be the
field of fractions. Replacing $F$ by its perfect closure in $Q$
and shrinking $S$, we may assume that $S$ is geometrically reduced.
Clearly, we may assume that our separable Hilbert subsets
$H_i\subset S$ are given by \'etale morphisms $T_i\ra S$.

Next, we reduce the problem to  a special case. Choose
a separating transcendence basis
$t_1,\ldots,t_n\in Q$ over $F$. Shrinking $S$, we obtain a separable
morphism
$f:S\ra\AA^n_F$.
 Choose 
nonempty open subsets 
$T_i'\subset T_i$  so that the composition
$T_i'\ra\AA^n_F$ is
\'etale, and let
$H_i'\subset\AA^n_F$ be the corresponding separable Hilbert subsets.
Then $f^{-1}(H_i')\subset H_i$, and we see that it suffices to treat
the special case
$S=\AA_F^n$.

We now proceed by induction on $n\geq 1$.
First, suppose $S=\AA_F^1$.
According to \cite{Fried; Jarden 1986},  Theorem 12.9,
there are nonempty open subsets $U_i\subset\AA^2_k$ with
$$
\left\{a+tb\mid a,b\in U_i(k)\right\}\subset H_i.
$$
Since $k$ is uncountable, the countable intersection
$\bigcap U_i$  contains uncountable many rational points
(choose a line not contained in any $\AA^2_k-U_i$).
Consequently $\bigcap H_i$ contains uncountable many rational points
as well.

Now suppose the result is already true for $n\geq 1$.
Let $F'$ be the function field of $\AA^n_F$, and consider
the projection $\AA^{n+1}_F\ra\AA^n_F$. Its generic fiber is
isomorphic to the affine line
$\AA^1_{F'}$. Set $H=\bigcap H_i$. We just saw that the
intersection $H\cap
\AA^1_{F'}$ contains a rational point
$z\in\AA^1_{F'}$. Its closure $S'\subset\AA^{n+1}_F$ is
$n$-dimensional. Applying our preliminary reduction and the 
induction hypothesis, we conclude that $H=\bigcap H_i$ contains
uncountably many closed points.
\qed

\begin{remark}
\mylabel{separable}
Under the additional hypothesis that $S$ is geometrically reduced,
the preceding proof   shows that there are uncountably many
closed points $s\in\bigcap H_i$ whose residue field extension
$F\subset \k(s)$ are separable.
\end{remark}

Summing up, we obtain the following result about
families of N\'eron--Severi groups:

\begin{theorem}
\mylabel{uncountable}
Let $k$ be an uncountable field, $k\subset F$ a nonalgebraic 
finitely generated field extension,
$S$ an integral $F$-scheme of finite type of dimension at least one,
and
$f:X\ra S$ a proper morphism. Then there are uncountably many closed
points
$s\in S$ with
$\rho(X_s)=\rho(X_\eta)$.
\end{theorem}

\proof
Combine Proposition \ref{countable} and Lemma \ref{classical
Hilbert}.
\qed

\medskip
The situation simplifies   if we look at Picard numbers
of   geometric fibers $X_{\bar{s}}=X\otimes\bar{\k}(s)$ instead
of   schematic fibers $X_s$.
More precisely:

\begin{theorem}
\mylabel{geometric Picard numbers}
Let $S$ be an integral noetherian scheme and
$f:X\ra S$ a proper morphism. Then there are countably many
nonempty open subsets $U_i\subset S$ such that
$\rho(X_{\bar{s}})= \rho(X_{\bar{\eta}})$ for all
$s\in\bigcap U_i$.
\end{theorem}

\proof
As in the proof for Proposition \ref{countable}, we reduce
to the case that $\Pic_{X/S}$ and $J\subset\Pic_{X/S}$ are
representable by group schemes, and set $\NS_{X/S}=\Pic_{X/S}/J$. For
each irreducible component $A_i\subset\NS_{X/S}$, choose a nonempty
open subset
$U_i\subset S$ over which $A_i$ is finite and flat. Then $\bigcap
U_i$ is the desired countable intersection of open subsets.

To see this, fix a point  $s\in \bigcap U_i$. Let $R=\O_{S,s}$ be
the corresponding local ring, $F=\k(s)$ its residue field, 
$F\subset F^\sep$ a separable closure,   and 
$R\subset R^\sh$ be the corresponding strict henselization. Then
$A_i\otimes R$ is
\'etale over
$R$, and $A_i\otimes R^\sh$ decomposes into disjoint sections by 
\cite{EGA IVd}, Proposition 18.5.19.
We conclude that $\NS(X\otimes F^\sep/F^\sep)$ is canonically
isomorphic to 
$\NS(X\otimes Q/Q)$,
where $Q$ is the function field of $R^\sh$.
The canonical mappings $\NS(X\otimes
F^\sep/F^\sep)\ra\NS(X\otimes\bar{F})$ and
$\NS(X\otimes Q/Q)\ra\NS(X\otimes\bar{Q})$ are obviously bijective,
so $\rho(X_{\bar{s}})=\rho(X_{\bar{\eta}})$.
\qed

\section{Fibered surfaces with small Picard number}
\mylabel{Fibered surfaces with small Picard number}

Fix an uncountable algebraically closed ground field $k$.
The task now is to construct fibered surfaces
with small Picard number defined over small transcendental extension
fields
$k\subset F$. The Tate--Shioda
formula then implies that their generic fibers have small Picard
group as well. The upshot is that there is a dense set of points
$x\in M_{g,n}$ with $\dim\overline{\left\{x\right\}}\leq 2$ such that
the marked points and the canonical class of the corresponding curve
$C_x$ generate $\Pic(C_x)\otimes\QQ$.

We start with some notation. A \emph{surface}  over $F$
is a proper 2-dimensional $F$-scheme $Y$ with $\Gamma(Y,\O_Y)=F$.
A \emph{fibration} on a surface $Y$ is curve $B$ together with a
proper morphism $f:Y\ra B$ satisfying $\O_B=f_*(\O_Y)$. The 
\emph{Tate--Shioda formula} relates the Picard number $\rho(Y)$ of a fibered
surface with the rank of the Mordell--Weil group of the
generic fiber (confer \cite{Tate 1966} and \cite{Shioda 1996}):

\begin{proposition}
\mylabel{tate shioda}
Let $Y$ be a normal surface whose singularities are $\QQ$-factorial,
and $f:Y\ra B$ be a fibration. The we have
$$
\rho(Y) - \rank\Pic^0(Y_\eta)/\Pic^0(Y) = 
2+\sum_{b\in B}(\rho(Y_b)-1).
$$
\end{proposition}

\proof
The intersection form is negative semidefinite on the group of Weil
divisors supported by a given fiber
$f^{-1}(b)\subset Y$, and the fiber $f^{-1}(b)$ generates the radical
over $\QQ$. Consequently, the vertical Weil divisors generate a subgroup
of rank $1+\sum_{b\in B}(\rho(Y_b)-1)$ inside the N\'eron--Severi group.

Next, choose a horizontal curve $H\subset Y$.
Given any Weil divisor $C$, we may subtract a suitable rational
multiple of $H$ until $C_\eta$ has degree zero.
Since $\Pic(Y)\ra \Pic(Y_\eta)$ has finite cokernel, we conclude that
the horizontal Weil divisors generate a subgroup of rank
$1+\rank\Pic^0(Y_\eta)/\Pic^0(Y)$ inside the N\'eron--Severi group.
This easily implies the formula.
\qed 

\medskip
This formula often allows us to control the rank of the Jacobian
$\Pic^0(Y_\eta)$ in terms of the Picard number $\rho(Y)$.

Now fix an integer $n\geq 0$ and a genus $g\geq 2$, and let $X$
be an $n$-pointed stable curve of genus $g$ over $k$.
For simplicity, we also assume $\Aut(X)=0$.
Then the closed point $x\in\qM_{g,n}$ corresponding to $X$ lies in
the smooth locus \cite{Lonstedt 1984} . As in the proof of Proposition
\ref{torsion free}, we blow up the center $x\in\qM_{g,n}$ and localize
at the generic point of the resulting exceptional divisor. This produces
a discrete valuation ring $A$ whose residue field $F=A/\maxid_A$ is a
purely transcendental field extension of    degree $3g-4+n$, and whose
field of fractions $A\subset Q$ is the function field of the moduli
space $\qM_{g,n}$.

Write $F=k(t_i)$ for some transcendence
basis
$t_i\in F$, and choose lifts
$t_i\in A$. We obtain an inclusion $k[t_i]\subset A$, and in turn a
lift
$F\subset A$. Then $F\subset Q$ is a finitely generated field
extension of transcendence degree one, which corresponds to a
proper normal curve
$B$ over $F$. The discrete valuation ring $A\subset Q$ defines a
rational point   $b_0\in B$, hence $F=\Gamma(B,\O_B)$, and $B$ is
geometrically integral.
The tautological curve $C\ra\qM_{g,n}$ defines a fibered
regular surface
$f:Y\ra B$. By construction, we have $f^{-1}(b_0)= X\otimes F$ and
$f^{-1}(\eta)= C_\eta$.
To apply the Tate--Shioda formula, we have to ensure the following:

\begin{lemma}
\mylabel{vanishing map}
Suppose the normalization $\tilde{X}$ of $X$ has only rational
components. Then the canonical map $\Pic^0(Y)\ra\Pic^0(Y_\eta)$
vanishes.
\end{lemma}

\proof
Let $F\subset\bar{F}$ be an algebraic closure. The normalization
$\tilde{Y}$ of $\bar{Y}=Y\otimes\bar{F}$ is geometrically normal,
so $\Pic^0_{\tilde{Y}/\bar{F}}$ is proper by
\cite{FGA VI}, Theorem 2.1. 
Moreover, it contains a  unique
abelian subscheme  with the same underlying topological space
by \cite{FGA VI}, Corollary 3.2. Its dual abelian variety is the
Albanese variety
$\Alb_{\tilde{Y}/\bar{F}}$. The resulting morphism
$\tilde{Y}\ra\Alb_{\tilde{Y}/\bar{F}}$ factors over $\tilde{B}$,
because the normalization of  $X$ has genus zero.

According to \cite{SGA 6}, Expos\'e XII, Corollary 1.5,
the map $\Pic^0_{\bar{Y}/\bar{F}}\ra \Pic^0_{\tilde{Y}/\bar{F}}$
is affine. It follows that  the cokernel of
$\Pic^0_{B/F}\ra\Pic^0_{Y/F}$ is affine.
Using that  $\Pic^0_{Y_\eta/\eta}$ is abelian, we infer that the map
$\Pic^0_{Y/F}\ra\Pic^0_{Y_\eta/\eta}$ factors over the origin, 
hence $\Pic^0(Y)\ra\Pic^0(Y_\eta)$ vanishes.
\qed

\medskip
From now on, we shall assume that the normalization of $X$
has only rational components.
Then Theorem \ref{Franchetta}, together  with the Tate--Shioda
formula gives
$$
\rho(Y) - (n+1) = 2+\sum_{b\in B}(\rho(Y_b)-1).
$$
We now use this formula to specialize $Y$ and maintain control over
$\Pic(Y_\eta)$.
Recall that $t_i\in F$, $1\leq i\leq 3g-4+n$ is a transcendence basis
over
$k$. Consider the rational function field $L=k(t_1)$. Then we may
view 
$L\subset F$ as the function field of $\PP_L^{3g-5+n}$. Our
$F$-schemes
$X$ and
$B$ extend to proper flat morphisms $\foY\ra S$ and $\foB\ra S$,
respectively, over some open subsets $S\subset \PP_L^{3g-5+n}$.
Shrinking $S$, we also have a morphism $f:\foY\ra\foB$ such that
$\O_\foB\ra f_*(\O_\foY)$ is bijective, and remains bijective after
any base change. Let $\foY_{\eta_s}\subset\foY_s$ be the generic
fiber of the induced fibration $\foY_s\ra\foB_s$.

\begin{proposition}
\mylabel{specialization}
There is an uncountable dense subset of closed
points $s\in S$ such that $\foY_{\eta_s}$ is a stable
$n$-pointed curve of genus
$g$, and that the marked points and the canonical class generate
$\Pic(\foY_{\eta_s})\otimes \QQ$.
\end{proposition}

\proof
Replacing $S$ by some nonempty open subset, we may assume that
the following holds: First, the fibers of $\foY\ra S$ and
$\foB\ra S$ are geometrically integral
by \cite{EGA IVc}, Theorem 9.7.7, and geometrically unibranch
by Proposition \ref{unibranch constructible}.
Second, the maps $\Pic^0_{\foY_s/s}\ra\Pic^0_{\foY_b/b}$ factor over
the origin for almost all closed points
$b\in\foB_s$. It then follows that
$\Pic^0(\foY_s)\ra\Pic^0(\foY_{\eta_s})$ vanishes. 
Third, the generic fibers $\foY_{\eta_s}$ are $n$-pointed smooth
curves of genus
$g$.

By Theorem \ref{uncountable}, there are uncountably many closed
points $s\in S$ with $\rho(\foY_s)=\rho(Y)$.
Let $\tilde{\foY}_s\ra\foY_s$ be the normalization and
$\bar{\foY}_s\ra\tilde{\foY}_s$ a resolution of singularities. Then
$\bar{\foY}_s$ is a fibered surface with generic
fiber $\bar{\foY}_{\eta_s}=\foY_{\eta_s}$.
We have $\rho(\tilde{\foY}_s)=\rho(\foY_s)$ by Proposition
\ref{finite}. The Picard number $\rho(\bar{\foY}_s)$ of the
desingularization is usually larger. However, any additionally
classes lie in the fibers of
$\bar{\foY}_s\ra\foB_s$.
The Tate--Shioda formula therefore implies that
$\Pic(Y_\eta)$ and $\Pic(\foY_{\eta_s})$ have the same rank.
Recall that $Y_\eta$ is the generic $n$-pointed curve of genus $g$.
Using  Theorem \ref{Franchetta}, 
we infer that the canonical class and the marked points generate
$\Pic(\foY_{\eta_s})\otimes \QQ$.
\qed

\medskip
We come to the second main result of this paper:

\begin{theorem}
\mylabel{dimension two}
Let $k$ be an uncountable algebraically closed field.
Then there is an uncountable dense   set of points
$x\in M_{g,n}$ with $\dim\overline{\left\{x\right\}}\leq 2$ such that
the canonical class and the marked sections generate
$\Pic(C_x)\otimes\QQ$, where $C_x$ is the curve corresponding to the
point $x\in M_{g,n}$.
\end{theorem}

\proof
Consider the family $\foY\ra\foB$ of fibered surfaces over $S$
constructed in Proposition \ref{specialization}.
We then have a rational map $\foB\dashrightarrow M_{g,n}\otimes L$
whose image is a divisor, because the composition
$\foB\dashrightarrow M_{g,n}$ is dominant.  After shrinking $S$, we
find an open subset $\foU\subset\foB$ whose fibers over $S$ are
nonempty, such that $\foU\ra M_{g,n}\otimes L$ is everywhere defined
and quasifinite.   This is because $\dim(\foB)=\dim(M_{g,n})-1$.

Consider the
uncountable dense set of points $s\in S$ from Proposition
\ref{specialization}, such that
$\Pic(\foX_{\eta_s})\otimes\QQ$ is generated by the canonical class
and the marked sections. Let $x\in M_{g,n}$ be the images of the
points
$\eta_s\in \foB$. By construction, the residue fields $\k(x)$ have
transcendence degree 
$\leq 2$ over $k$. Hence the $x\in M_{g,n}$ constitute
the desired uncountable dense   subset.
\qed
\section{Open problems}

We close the paper by listing some open problems:

\medskip
(1) Theorem \ref{dimension two} states that the marked points
and the canonical class generate $\Pic(C_x)$. If $l\in\Pic(C_x/x)$ is
a rational point, what denominators are necessary to write $l$ as a
linear combination of the canonical class and the marked points?
Is it sufficient to allow powers of 
the characteristic exponent $p\geq 1$ as denominators?

\medskip
(2) Is it possible to choose a dense set of nonclosed points 
$x\in M_{g,n}$ as in Theorem
\ref{dimension two} so that all closures
$\overline{\left\{ x \right\}}$ are 1-dimensional?
This seems to rely on improved versions of Hilbert's Irreducibility
Theorem for prime fields.

\medskip
(3) Can we say more about the structure of countable intersections
of nonempty Hilbert sets in Theorem \ref{uncountable}?
Obviously, we cannot expect such sets to have a reasonable
algebraic structure. Is it possible to write such
sets as a countable union of sets with some
sort of algebraic structure?

\medskip
(4) Does the Strong Franchetta Conjecture generalize to
other moduli problems, and if so, in what form?
What about polarized abelian varieties, or surfaces of general type,
or canonically polarized varieties of higher dimensions?
For example, Silverberg \cite{Silverberg 1985} 
showed that generic complex abelian
varieties have finite Mordell--Weil groups.



\begin{thebibliography}{ccccc}

\bibitem{Arbarello et al 1985}
E.\ Arbarello, M.\ Cornalba, P.\ Griffiths, J.\ Harris:
Geometry of algebraic curves. I. 
Grundlehren Math.\ Wiss.\  267. 
Springer, New York, 1985.

\bibitem{Arbarello; Cornalba 1987}
E.\ Arbarello, M.\ Cornalba:
The Picard groups of the moduli spaces of curves. 
Topology 26 (1987), 153--171.

\bibitem{Arbarello; Cornalba 1998}
E.\ Arbarello, M.\ Cornalba:
Calculating cohomology groups of moduli spaces of curves via 
algebraic geometry. 
Inst.\ Hautes Etudes Sci.\ Publ.\ Math.\  88 (1998), 97--127.

\bibitem{Belyi 1979}
G.\ Belyi:
Galois extensions of a maximal cyclotomic field. 
Izv.\ Akad.\ Nauk SSSR Ser.\ Mat.\ 43 (1979),  267--276.

\bibitem{Boggi; Pikaart 2000}
M.\ Boggi, M.\ Pikaart: 
Galois covers of moduli of curves.
Compositio Math.\ 120 (2000),  171--191.

\bibitem{Bosch; Luetkebohmert; Raynaud 1990}
S.~Bosch, W.~L\"utkebohmert, M.~Raynaud:
N\'eron models.
Ergeb.\ Math.\ Grenzgebiete (3) 21.
Springer, Berlin, 1990.

\bibitem{A 4-7}
N.~Bourbaki:
Algebra II. Chapters 4--7.
Springer, Berlin, 1990.

\bibitem{Deligne; Mumford 1969}
P.\ Deligne, D.\ Mumford:
The irreducibility of the space of curves of given genus. 
Publ.\ Math.\, Inst.\ Hautes \'Etud.\ Sci.\  36 (1969), 75--109.

\bibitem{Deligne 1985}
P.\ Deligne:
Le lemme de Gabber. 
Ast\'erisque  127 (1985), 131--150. 

\bibitem{Franchetta 1954}
A.\ Franchetta:
Sulle serie lineari razionalmente determinate sulla 
curva a moduli generali di dato genere. 
Matematiche (Catania) 9 (1954), 126--147.

\bibitem{Fried; Jarden 1986}
M.\ Fried, M.\ Jarden: 
Field arithmetic. 
Ergeb.\ Math.\ Grenzgebiete (3) 11. 
Springer, Berlin, 1986.

\bibitem{Fried; Klassen; Kopeliovich 2001}
M.\ Fried, E.\ Klassen, Y.\ Kopeliovich:
Realizing alternating groups as monodromy groups of genus one covers.
Proc.\ Amer.\ Math.\ Soc.\ 129 (2001),  111--119.

\bibitem{Fulton 1969}
W.\ Fulton:
Hurwitz schemes and irreducibility of moduli of algebraic curves.
Ann.\ Math.\ 90 (1969), 542--575.

\bibitem{EGA I}
A.\ Grothendieck, J.A.\ Dieudonn\'e:
\'El\'ements de g\'eom\'etrie alg\'ebrique I:
Le langage de sch\'emas.
Grundlehren Math.\ Wiss.\ 166.
Springer, Berlin, 1970.

\bibitem{EGA IIIa}
A.\ Grothendieck:
\'El\'ements de g\'eom\'etrie alg\'ebrique III:
\'Etude cohomologique des faiscaux coh\'erent.
Publ.\ Math., Inst.\ Hautes \'Etud.\ Sci.\ 11 (1961).

\bibitem{EGA IVb}
A.\ Grothendieck:
\'El\'ements de g\'eom\'etrie alg\'ebrique IV: \'Etude locale des
sch\'emas et de morphismes de sch\'emas.
Publ.\ Math., Inst.\ Hautes \'Etud.\ Sci.\   24 (1965).

\bibitem{EGA IVc}
A.\ Grothendieck:
\'El\'ements de g\'eom\'etrie alg\'ebrique IV: \'Etude locale des
sch\'emas et de morphismes de sch\'emas.
Publ.\ Math., Inst.\ Hautes \'Etud.\ Sci.\  28 (1966).

\bibitem{EGA IVd}
A.\ Grothendieck:
\'El\'ements de g\'eom\'etrie alg\'ebrique IV: \'Etude locale des
sch\'emas et de morphismes de sch\'emas.
Publ.\ Math., Inst.\ Hautes \'Etud.\ Sci.\   32 (1967).

\bibitem{FGA V}
A. Grothendieck:
Les sch\'emas de Picard: th\'eor\'emes d'existence.
S\'eminaire Bourbaki, Exp.\ 232 (1962), 143--161.

\bibitem{FGA VI}
A.\ Grothendieck: 
Les sch\'emas de Picard:
propri\'et\'es g\'en\'erales. 
S\'eminaire Bourbaki, Exp.\ 236 (1962), 221--243. 

\bibitem{GB II}
A.\ Grothendieck:
Le groupe de Brauer II. 
In: J.~Giraud (ed) et al.: Dix exposes sur la cohomologie des schemas, pp.\
88--189.
North-Holland, Amsterdam, 1968. 

\bibitem{GB III}
A.~Grothendieck:
Le groupe de Brauer III.
In: J.~Giraud (ed.) et al., Dix exposes sur la cohomologie des schemas, pp.\
88--189.
North-Holland, Amsterdam, 1968.

\bibitem{SGA 1}
A.~Grothendieck et al.:
Rev\^etements \'etales et groupe fondamental.
Lect.\ Notes Math.\  224.
Springer, Berlin, 1971.

\bibitem{SGA 3a}
A.~Grothendieck et al.:
Schemas en groupes I.
Lect.\ Notes Math.\  151.
Springer, Berlin, 1970.

\bibitem{SGA 6}
A.\ Grothendieck et al.:
Th\'eorie des intersections et th\'eor\`eme de Riemann-Roch.
Lect.\ Notes   Math.\ 225.
Springer, Berlin, 1971.

\bibitem{Harer 1983}
J.\ Harer:
The second homology group of the mapping class group of an orientable surface. 
Invent.\ Math.\ 72 (1983), 221--239.

\bibitem{Hartshorne 1977}
R.~Hartshorne:
Algebraic geometry.
Grad.\ Texts   Math.\ 52.
Springer, Berlin,  1977.

\bibitem{Huebl; Kunz 1994}
R.\ H\"ubl, R.\ Kunz:
On algebraic varieties over fields of prime characteristic. 
Arch.\ Math.\ (Basel) 62 (1994), 88--96. 

\bibitem{Jouanolou 1983}
J.-P.\ Jouanolou:
Th\'eor\`emes de Bertini et applications.
Progress in Mathematics 42. 
Birkh\"auser, Boston, 1983.

\bibitem{Kleiman 1966}
S.~Kleiman:
Toward a numerical theory of ampleness.
Ann.\ Math.\  84  (1966), 293--344.

\bibitem{Kodaira 1963}
K.\ Kodaira:
On compact analytic surfaces II.
Ann.\ of Math.\  77 (1963), 563--626.

\bibitem{Kollar 1997}
J.\ Koll\'ar:
Quotient spaces modulo algebraic groups. 
Ann.\ of Math.\ 145 (1997), 33--79.

\bibitem{Kouvidakis 1991}
A.\ Kouvidakis: 
The Picard group of the universal Picard varieties over the 
moduli space of curves. 
J.\ Differential Geom.\ 34 (1991), 839--850. 

\bibitem{Lang 1983}
S.~Lang:
Fundamentals of diophantine geometry.
Springer, Berlin, 1983.

\bibitem{Lonstedt 1984}
K.\ L\o nsted:
The singular points on the moduli spaces for smooth curves. 
Math.\ Ann.\ 266 (1984), 397--402. 

\bibitem{Loojienga 1994}
E.\ Loojienga:
Smooth Deligne--Mumford compactifications by meands of Prym level
structures.
J.\ Alg.\ Geom.\ 3 (1994), 283--293.

\bibitem{Mestrano 1987}
N.\ Mestrano:
Conjecture de Franchetta forte. 
Invent.\ Math.\ 87 (1987), 365--376.

\bibitem{Milne 1980}
J.\ Milne: \'Etale cohomology. 
Princeton Mathematical Series, 33. 
Princeton University Press, Princeton, 1980. 

\bibitem{Moriwaki 2001}
A.\ Moriwaki:
The $\QQ$-Picard group of the moduli space of curves in positive
characteristic. 
Internat.\ J.\ Math.\ 12 (2001), 519--534.

\bibitem{Morrison 1985}
D.\ Morrison:
The birational geometry of surfaces with rational double points. 
Math.\ Ann.\ 271 (1985), 415--438.

\bibitem{Mumford; Fogarty; Kirwan 1993}
D.~Mumford, J.~Fogarty, F.~Kirwan:
Geometric invariant theory. Third edition.
Ergeb.\ Math.\  Grenzgebiete (3) 34. Springer, Berlin, 1993.

\bibitem{Pikaart; de Jong 1994}
M.\ Pikaart, A.\ de Jong:
Moduli of curves with non-abelian level structure.
In: R.\ Dijkgraaf, C.\ Faber, G.\ van der Geer (eds.),
The moduli space of curves, pp.\ 483--509.
Progress in Math.\ 129.
Birkh\"auser, Boston, 1995.

\bibitem{Poonen 2000}
B.\ Poonen:
Varieties without extra automorphisms I. Curves. 
Math.\ Res.\ Lett.\ 7 (2000),  67--76.

\bibitem{Reid 1994}
M.\ Reid:
Nonnormal del Pezzo surfaces.
Publ.\ Res.\ Inst.\ Math.\ Sci.\ 30 (1994), 695--727.

\bibitem{Saidi 1997}
M.\ Sa{\"\i}di:
Rev\^etements mod\'er\'es et groupe fondamental de graphe de groupes.
Compositio Math.\ 107 (1997), 319--338.

\bibitem{Serre 1972}
J.-P.\ Serre:
Cohomologie galoisienne.
Fifth edition. Lect.\ Notes  Math.\ 5. 
Springer, Berlin, 1994.

\bibitem{Serre 1989}
J.-P.\ Serre:
Lectures on the Mordell-Weil theorem.
Vieweg, Braunschweig, 1989.

\bibitem{Shioda 1996}
T.\ Shioda:
Mordell--Weil lattices for higher genus fibration
over a curve. 
In: K.\ Hulek (ed.) et al.,
New trends in algebraic geometry, pp. 359--373.
London Math.\ Soc.\ Lecture Note Ser.\ 264.
Cambridge Univ.\ Press, Cambridge, 1999. 

\bibitem{Silverberg 1985}
A.\ Silverberg:
Mordell--Weil groups of generic abelian varieties. 
Invent.\ Math.\ 81 (1985), 71--106. 

\bibitem{Tate 1966}
J.\ Tate:
On the conjectures of Birch and Swinnerton-Dyer 
and a geometric analog. 
S\'eminaire Bourbaki, Exp.\ 306 (1966).

\end{thebibliography}
\end{document}